\documentclass[11pt]{article}
\usepackage[title]{appendix}
\usepackage{xcolor}
 \usepackage{geometry}                % See geometry.pdf to learn the layout options. There are lots.
\usepackage{graphicx}
\usepackage{subcaption}
\usepackage{pdflscape}
\usepackage{amssymb}
\usepackage[normalem]{ulem}
\usepackage{hyperref}
\usepackage{enumitem}
\usepackage{mathrsfs}
\newcounter{casecounter}
\usepackage{epstopdf}
\usepackage{rotating}
\usepackage{longtable} % for 'longtable' environment
\usepackage{adjustbox}
\usepackage{float}

\usepackage{color}
\usepackage{bbm, dsfont}
\usepackage{pst-node}
\usepackage{tikz-cd}
\usepackage{amsfonts} %mathbb{}
\usepackage{geometry}
\usepackage{amsthm}
\usepackage{amsmath}
\usepackage{titlesec}
\usepackage{fixltx2e}
\usepackage{amssymb}
\usepackage{enumitem}
\usepackage{float}
\restylefloat{table}
\usepackage{tikz}
\usetikzlibrary{tikzmark}
\usetikzlibrary{arrows}
\usetikzlibrary{cd}
\usetikzlibrary{automata, positioning, arrows}
\usepackage [english]{babel}
\usepackage [autostyle, english = american]{csquotes}
\usepackage{algorithm}
\usepackage[noend]{algpseudocode} %pseudocode
\makeatletter
\def\BState{\State\hskip-\ALG@thistlm}
\makeatother

\def\X^n{\text{X}^n}

\usepackage{hyperref}
\usepackage[normalem]{ulem}
\newlist{casess}{enumerate}{1}
\setlist[casess]{label=     \textbf{Case} \arabic*:}
\usepackage{mathtools}

\DeclarePairedDelimiter\floor{\lfloor}{\rfloor}
\newcommand\addtag{\refstepcounter{equation}\tag{\theequation}}
\makeatletter
\newcommand*{\rom}[1]{\expandafter\@slowromancap\romannumeral #1@}
\makeatother

\theoremstyle{theorem}
\newtheorem{thm}{Theorem}[section]
\newtheorem{prop}[thm]{Proposition}
\newtheorem{cor}[thm]{Corollary}
\newtheorem{lemma}[thm]{Lemma}

\theoremstyle{definition}
\newtheorem{defn}{Definition}[section]
\newtheorem{rmk}[thm]{Remark}

\newtheorem{exmp}[thm]{Example}

\usepackage{etoolbox}

\makeatletter
\patchcmd{\ttlh@hang}{\parindent\z@}{\parindent\z@\leavevmode}{}{}
\patchcmd{\ttlh@hang}{\noindent}{}{}{}
\makeatother

\usepackage{listings}
\usepackage{color} %red, green, blue, yellow, cyan, magenta, black, white
\definecolor{mygreen}{RGB}{28,172,0} % color values Red, Green, Blue
\definecolor{mylilas}{RGB}{170,55,241}

\newlist{Assumptions}{enumerate}{1}
\setlist[Assumptions]{label=     \textbf{Assumption} \arabic*:}

\makeatletter

\newsavebox{\@brx}
\newcommand{\llangle}[1][]{\savebox{\@brx}{\(\m@th{#1\langle}\)}%
  \mathopen{\copy\@brx\kern-0.5\wd\@brx\usebox{\@brx}}}
\newcommand{\rrangle}[1][]{\savebox{\@brx}{\(\m@th{#1\rangle}\)}%
  \mathclose{\copy\@brx\kern-0.5\wd\@brx\usebox{\@brx}}}
\makeatother

\usepackage{lipsum} % for generating filler text
\usepackage{titlesec}
\titleformat{\subsection}[runin]% runin puts it in the same paragraph
       {\normalfont\bfseries}% formatting commands to apply to the whole heading
       {\thesubsection}% the label and number
       {0.5em}% space between label/number and subsection title
       {}% formatting commands applied just to subsection title
       [.]% punctuation or other commands following subsection title

\begin{document}
\title{Growth of finitely generated subgroups of the topological full groups of inverse semigroups of bounded type}
\author{Zheng Kuang\thanks{Texas A\&M University, College Station, USA; Email: kzkzkzz@tamu.edu}}
\date{}
\maketitle

\begin{abstract}
 Given an inverse semigroup $G_0$ of bounded type, we show, along with some other assumptions, that if the set of incompressible elements of $G_0$ is finite, then any finitely generated subgroup $G$ of the topological full group $\mathsf{F}(G_0)$ that is orbit equivalent to $G_0$ has subexponential growth with a bounded power in the exponent. 
\end{abstract}
\tableofcontents
\section{Introduction}
In the previous work \cite{kua24a}, we explored a class of groups of intermediate growth, namely, groups with finitely many incompressible elements, applying the methods of \cite{BNZ} to groups whose orbital graphs are not quasi-isometric to a line. The main result established in it is that for a group $G$ of bounded type, acting minimally on the space of infinite paths of a simple and stationary Bratteli diagram and defined by a stationary automaton, if the set of incompressible elements of 
$G$ is finite, then the growth function of $G$, denoted $\gamma_G(R)$, is bounded above by a function of the form $\exp(CR^\alpha)$ for $\alpha\in (0, 1)$. 

In this paper, we continue the study of the growth of groups with finitely many incompressible elements. We generalize the main result from \cite{kua24a} in the following aspects. 

First, we drop the assumption that the Bratteli diagram $\mathsf{B}$ and the automaton are stationary. We consider the tile inflation processes on Bratteli diagrams $\mathsf{B}$ defining inverse semigroups $G_0$ of bounded type. We will consider tile inflation satisfying \emph{expansion} condition, which ensures polynomial growth of the orbital graphs. The latter is essential for obtaining subexponential growth estimate. The expansion is also one of the key conditions ensuring \emph{linear repetitivity} of tiles.  We will show that if the tile inflation process is expanding, then one of the necessary conditions for linearly repetitivity is that the number of edges on each level of $\mathsf{B}$ is uniformly bounded. This exhibits a connection with \textit{continued fractions of bounded type}, as we will see in Section \ref{dihegrp}. 
  
 Next, we present a systematic method of finding groups of bounded type from an inverse semigroup $G_0$ of bounded type. In \cite{kua24a}, examples of groups of bounded type were only obtained when the tile inflation processes directly define a group or the inverse semigroup $G_0$ can be directly turned into a group $G$. See the example in \cite[Subsection 5.3]{kua24a}. In this paper, we show that if $G_0$ is an inverse semigroup of bounded type, then any finitely generated subgroup $G$ of the topological full group $\mathsf{F}(G_0)$ that is orbit equivalent to $G_0$ is a group of bounded type. If $G_0$ has finitely many incompressible elements, then so does $G$, which implies that $G$ is periodic (torsion). Periodicity and the main theorem together imply that the group $G$ has intermediate growth. 
 
 The main theorem of this paper is as follows (Theorem \ref{main}).  

 \begin{thm}
   Let $\mathsf{B}$ be a simple Bratteli diagram. Suppose the inverse semigroup $G_0=\langle \mathcal{S}\rangle$ is defined by an expanding and linearly repetitive tile inflation process and acts minimally on $\Omega(\mathsf{B})$. Let $G<\mathsf{F}(G_0)$ be a finitely generated subgroup that is orbit equivalent to $G_0$. If the set of incompressible elements of $G_0$ is finite, then there exists $\alpha\in (0,1)$ such that the growth function $\gamma_G$ of $G$ satisfies $\gamma_G(R)\preccurlyeq \exp(R^{\alpha})$. 
 \end{thm}

 The paper is organized as follows. Section \ref{prelim} reviews basic notions that will be used in this paper. Section \ref{exptile} discusses expanding tile inflations, linear repetitivity of tiles, and polynomial growth of orbits. We show that the condition that a tile inflation process is expanding implies that the diameters of finite tiles grow exponentially and thus infinite tiles have polynomial growth and finite tiles on the same levels of $\mathsf{B}$ have asymptotically equivalent diameters. The last property is used to show linear repetitivity. Section \ref{prvmain} is the proof of the main theorem. The method of counting ``traverses" developed by L.~Bartholdi, V.~Nekrashevych and T.~Zheng \cite{BNZ} is the main tool in the proof. In Section \ref{dihegrp}, we describe a family of infinite dihedral groups $D_{\infty}$ to exhibit the connection between groups of bounded type and continued fractions of bounded type. In Section \ref{nonlin13}, we provide a family of groups with non-linear orbital graphs. We fragment some of the groups in Sections \ref{dihegrp}, \ref{nonlin13} to produce the explicit examples of groups of intermediate growth acting on the spaces of infinite paths of non-stationary Bratteli diagrams.

\section{Preliminaries}\label{prelim}
Throughout this paper, we use right action notation. The notation $G_0$ is reserved for inverse semigroups. A product $g=F_1F_2\ldots F_m\in G_0$ is read from left to right. Paths in the Bratteli diagram $\mathsf{B}$, e.g. $\gamma=\ldots e_2e_1$, are read from right to left. Action of a group or semigroup element on a path is denoted $\gamma\cdot g$. All the tiles are assumed to be connected graphs. 

\subsection{Tile inflations and automata}\label{tia}
\begin{defn}\label{graphs}
A \textit{graph with boundary} $\Gamma$ consists of a set of \textit{vertices} $V(\Gamma)$, a set of \textit{edges} $E(\Gamma)$,  two partially defined maps
\[\mathsf{s}:E(\Gamma)\dashrightarrow V(\Gamma),\] 
\[\mathsf{r}:E(\Gamma)\dashrightarrow V(\Gamma),\]
and a map
\[E(\Gamma)\rightarrow E(\Gamma),\text{    } e\mapsto e^{-1},\] which satisfy the conditions: 
\begin{enumerate}
\item each $e\in E(\Gamma)$ belongs to at least one of the domains of $\mathsf{s}$ and $\mathsf{r}$; 
\item for each $e\in E(\Gamma)$ we have $(e^{-1})^{-1}=e$, $e^{-1}\neq e$ and $\mathsf{s}(e)=\mathsf{r}(e^{-1})$ or $\mathsf{r}(e)=\mathsf{s}(e^{-1})$. 
\end{enumerate}
The maps $\mathsf{s},\mathsf{r}$ are called respectively \textit{source} and \textit{range} maps, which are defined on subsets of $E(\Gamma)$. For each $e\in E(\Gamma)$, the edge $e^{-1}$ is called the \textit{inverse} edge of $e$. 
\end{defn}

We interpret an edge $e$ as an arrow from $v_0=\mathsf{s}(e)$ to $v_1=\mathsf{r}(e)$ and sometimes denote them as $v_0\rightarrow v_1$.

Informally, we allow some of the edges to be ``hanging'' with source or range to be not in the set of vertices.

Let $\Gamma$ be a graph with boundary. Let us define a class of edge-labeling.

Let $\mathcal{S}$ be a finite set of labels. We assume that $\mathcal{S}$ is symmetric, i.e., there is an involution map \[\mathcal{S}\rightarrow\mathcal{S}, \textup{ } F\mapsto F^{-1}\] mapping each $F\in\mathcal{S}$ to a unique $F^{-1}\in\mathcal{S}$ with $(F^{-1})^{-1}=F$. 

Each element in the set of edges $E(\Gamma)$ is labeled by an element of $\mathcal{S}$, where if $e\in E(\Gamma)$ is labeled by $F\in\mathcal{S}$, then $e^{-1}$ is labeled by the unique $F^{-1}\in\mathcal{S}$ such that one of the following cases will happen.
\begin{enumerate}
    \item There exist $\gamma_1,\gamma_2\in V(\Gamma)$ connected by an arrow (edge) $\gamma_1\rightarrow\gamma_2$ labeled by $F$, and by an arrow $\gamma_1\leftarrow\gamma_2$ labeled by $F^{-1}$ while no other edges connecting $\gamma_1,\gamma_2$ (if exist) are labeled by $F$ or $F^{-1}$. 
    \item\label{bdd} There exists $\gamma\in V(\Gamma)$ and $e\in E(\Gamma)$ such that $\mathsf{s}(e)=\gamma$ but $e$ is not in the domain of $\mathsf{r}$ and $e$ is labeled by $F$. Correspondingly, $e^{-1}$, labeled by $F^{-1}$, is such that $\mathsf{r}(e^{-1})=\gamma$ and it is not in the domain of $\mathsf{s}$. In other words, there exists a $\gamma\in V(\Gamma)$
    to which an outgoing edge labeled by $F$ and an incoming edge labeled by $F^{-1}$ are attached, but $\gamma$ is not connected to other vertices by these two edges. 
\end{enumerate}
If both 1 and 2 happen for the same $F$, then $\gamma_1\neq \gamma$.

\begin{defn}\label{bddpt}
   The point $\gamma$ in Case \ref{bdd} above is called a \textit{boundary vertex} of a graph. The edges labeled by $F$ and $F^{-1}$ at $\gamma$ are called \textit{boundary edges}. 
\end{defn} 

\begin{defn}\label{welllabel}
    Let $\Gamma$ be a graph whose edges are labeled by elements of the set $\mathcal{S}$. A vertex $v$ of $\Gamma$ is called \textit{well-labeled} if, for every $F\in \mathcal{S}$, there is at most one edge starting in $v$ labeled by $F$ and at most one edge ending in $v$ labeled by $F$. The graph $\Gamma$ is said to be \textit{well-labeled} if each of its vertices is well-labeled. It is said to be \textit{perfectly-labeled} if for every vertex $v$ and every $F\in \mathcal{S}$ there is exactly one edge starting in $v$ labeled by $F$ and exactly one edge ending in $v$ labeled by $F$. 
\end{defn}

\begin{defn}\label{brat} A \textit{Bratteli diagram} $\mathsf{B}$ consists of sequences $(V_1,V_2,\ldots)$ and $(E_1,E_2,\ldots)$ of finite sets and sequences of maps $\mathbf{s}_n:E_n\rightarrow V_n$ and $\mathbf{r}_n:E_n\rightarrow V_{n+1}$. The sets $\bigsqcup_{n\geq 1}V_n$ and $\bigsqcup_{n\geq 1}E_n$ are called, respectively, set of \textit{vertices} and set of \textit{edges} of the diagram. The maps $\mathbf{s}_n$ and $\mathbf{r}_n$ are called \textit{source map} and \textit{range map}, respectively, and they are assumed to be surjective. We will write $\mathbf{s}$ and $\mathbf{r}$ if no ambiguity would be caused. Denote $\mathsf{B}=((V_n)_{n=1}^{\infty},(E_n)_{n=1}^{\infty},\mathbf{s},\mathbf{r})$. Also denote by $\Omega(\mathsf{B})$ the space of all infinite paths of $\mathsf{B}$ starting in $V_1$, by $\Omega_n(\mathsf{B})$ the space of all paths starting in $V_1$ ending in $V_{n+1}$, and by $\Omega^*(\mathsf{B}):=\bigcup_{n=1}^{\infty}\Omega_n(\mathsf{B})$ the space of all finite paths. Here a \emph{path} in the diagram is a sequence $e_n\ldots e_2e_1$ of edges such that $\mathsf{r}(e_i)=\mathsf{s}(e_{i+1})$.
\end{defn}
The space $\Omega(\mathsf{B})$ is a closed subset with the subset topology of the direct product $\prod\limits_{n=1}^{\infty}E_n$. It is compact, totally disconnected, and metrizable.

A \textit{tile} $\mathcal{T}_{v,n}$ is a graph with boundary defined by a \textbf{tile inflation} process, described below, where $V(\mathcal{T}_{v,n})$ consists of all paths of length $n$ of $\mathsf{B}$ ending in the same vertex $v\in E_{n+1}$.

Let us describe this inflation process.
%which is a generalization of \cite[Definition 23]{Bon11}.
Suppose  all tiles up to the $n$-th level of $\mathsf{B}$ are constructed. 

\begin{defn}\label{admissible}
  Let $\mathcal{T}_{1},\mathcal{T}_2$ be tiles (possibly $\mathcal{T}_1=\mathcal{T}_2$) on the $n$-th level of $\mathsf{B}$. A pair of boundary edges $e_1, e_2$ of $\mathcal{T}_1$, $\mathcal{T}_2$, respectively are said to be \textit{compatible} if  both edges are labeled by the same label $F$ and either $\mathsf{s}(e_1)$ and $\mathsf{r}(e_2)$ are undefined, or vice versa. The corresponding vertices of the edges $e_1$ and $e_2$ are called \emph{admissible}.
\end{defn}
\begin{exmp}
  If each $F=F^{-1}$, then each boundary point on a finite tile is admissible to itself. 
\end{exmp}

Let $v\in V_{n+2}$ and $\mathcal{T}_v$ be an $(n+1)$-st level tile. It is obtained as follows. 

We choose a set  $P_{n+1}$ (called the set of \textit{connectors} of level $n+1$) of triples 
$(e_1\gamma_1,e_2\gamma_2, F)$ where vertices $\gamma_1,\gamma_2$ of $n$th level tile are such that are admissible boundary points with corresponding compatible edges labeled by $F$ (the edges are unique, by our condition that the graphs are well labeled). The points $e_1\gamma_1$ and $e_2\gamma_2$ are called \textit{connecting points}. 

For each edge $e\in \mathbf{r}^{-1}_{n+1}(v)$, take a copy, denoted $\mathcal{T}_{e}$, of the tile $\mathcal{T}_{\mathbf{s}_{n+1}(e)}$. 

Consider a pair of edges $e_1,e_2\in\mathbf{r}_{n+1}^{-1}(v)$, and suppose that $\gamma_1$ and $\gamma_2$ are boundary points of $\mathcal{T}_{\mathbf{s}_{n+1}(e_1)}$ 
and $\mathcal{T}_{\mathbf{s}_{n+1}(e_2)}$ respectively, such that $(e_1\gamma_1,e_2\gamma_2, F)\in P_{n+1}$. 

Connect then the vertex $e_1\gamma_1$ of $\mathcal{T}_{e_1}$ to the vertex $e_2\gamma_2$ of $\mathcal{T}_{e_2}$ by a pair of arrows $e_1\gamma_1\rightarrow e_2\gamma_2$ and $e_1\gamma_1\leftarrow e_2\gamma_2$ labeled by $F$ and $F^{-1}$, respectively, or labeled by $F^{-1}$ and $F$, respectively. 

The tile $\mathcal{T}_v$ is obtained by applying the above process to all connectors in $P_{n+1}$. 
\begin{rmk}
    Note that a boundary point $\gamma$ on level $n$ (i.e., a path in the Bratteli diagram) must be continued to either a boundary point $e\gamma$ or a connecting point $e'\gamma$ on level $n+1$, or both happen. Meanwhile, non-boundary points cannot be continued to boundary points or connecting points. 
\end{rmk}
  Let $\gamma\in\Omega(\mathsf{B})$. Write $\gamma=...e_n\ldots e_2e_1$. Denote by $\gamma_n=e_n\ldots e_1$ its $n$-th \textit{truncation} ending in a vertex $v_{n+1}\in V_{n+1}$. 
  \begin{defn}\label{inftile}
      The \textit{infinite tile} of $\gamma$, denoted $\mathcal{T}_{\gamma}$, is the inductive limit of the embeddings of graphs
\begin{center}
    \begin{tikzcd}
\mathcal{T}_{v_2} \arrow[r, "\phi_{v_2,e_2}"] & \mathcal{T}_{v_3} \arrow [r, "\phi_{v_3,e_3}"] & \ldots \arrow [r, "\phi_{v_{n-1},e_{n-1}}"]& \mathcal{T}_{v_{n}} \arrow [r, "\phi_{v_n,e_{n}}"] & \ldots,
\end{tikzcd}
\end{center}
where $\phi_{v_i,e_{i}}:\eta\mapsto e_{i}\eta$, for $\eta\in \Omega_{i}(\mathsf{B})$ ending in $v_{i}$, is an ismorphism of $\mathcal{T}_{v_i}$ with a subgraph of $\mathcal{T}_{v_{i+1}}$. 
\end{defn}
\begin{defn}
    Let $\mathcal{T}$ be an infinite tile. An infinite path $\xi\in\mathcal{T}$ is said to be a \textit{boundary point of $\mathcal{T}$} if all its $n$-th truncations are boundary vertices in the finite tiles, for $n\in\mathbb{N}$.
\end{defn}

\begin{defn}
    Let $\textup{X}$ be a topological space. An \textit{inverse semigroup} acting on $\textup{X}$ by partial  homeomorphisms, denoted $\textup{X}\curvearrowleft G_0$, is a collection of homeomorphisms $g:\textup{Dom}(g)\rightarrow\textup{Ran}(g)$ between open subsets of $\textup{X}$ closed under composition and taking inverses.
\end{defn}

The labels of the edges of infinite tiles define bijections between subsets of $\Omega(\mathsf{B})$, as it follows from the construction of tile inflation and Definition \ref{welllabel}). Namely, for every label $F\in\mathcal{S}$ and an infinite path $\gamma\in\Omega(\mathsf{B})$ we define $\gamma\cdot F$ as $\mathsf{r}(e)$, where the edge $e$ of $\mathcal{T}_\gamma$ is such that $\mathsf{s}(e)=\gamma$ and $e$ is labeled by $F$. Note that such a (non-boundary) edge may not exist. On the other hand, if it exists, then it exists for all paths $\gamma'\in\Omega(\mathsf{B})$ that have a sufficiently long common beginning with $\gamma$. It follows that the domain of $F$ is open. It follows that every tile inflation rule defines a finite collection of partial homeomorphisms of $\Omega(\mathsf{B})$, hence generates an inverse semigroup.

On the other hand, the union of the domain of $F$ with the set of boundary points of infinite tiles that have a boundary edge labeled by $F$ is clopen.

In some cases the local homeomorphism of $\Omega(\mathsf{B})$ defined by a label $F$ can be extended to a homeomorphism of the closure of the domain. In other words, we can  connect the corresponding infinite tiles at the boundary points by the arrows labeled by $F$ and $F^{-1}$ to form a new well labeled graph. Then each $F\in\mathcal{S}$ defines a transformation (partial homeomorphism) with a \textbf{clopen domain} of $\Omega(\mathsf{B})$, denoted $F:\textup{Dom}(F)\rightarrow\textup{Ran}(F)$. Since each $F\in\mathcal{S}$ corresponds to a unique inverse $F^{-1}$, the set $\mathcal{S}$ actually generates an inverse semigroup, denoted $G_0=\langle \mathcal{S}\rangle$. 

\begin{defn}
Let $\gamma\in\Omega(\mathsf{B})$. The \textit{orbital graph} of $\gamma$, denoted $\Gamma_{\gamma}$, has the orbit $\gamma\cdot G_0$ as its set of vertices, and for every $\eta\in \gamma\cdot G_0$ and every $F\in \mathcal{S}$ defined at $\eta$, there is an arrow from $\eta$ to $\eta\cdot F$ labelled by $F$.
\end{defn}

\begin{rmk}
Let $\mathcal{T}_{\gamma}$ be an infinite tile. Then $\Gamma_{\gamma}=\mathcal{T}_{\gamma}$ if $\mathcal{T}_{\gamma}$ does not have any boundary points. If $\mathcal{T}_{\gamma}$ contains a boundary point $\xi$ corresponding to an $F\in\mathcal{S}$, and $F$ can be extended continuously at $\xi$, then we can ``paste" the tiles $\mathcal{T}_{\xi}(=\mathcal{T}_{\gamma})$ with $\mathcal{T}_{\xi\cdot F}$ by arrows $\xi\rightarrow \xi\cdot F$ labeled by $F$ and $\xi\leftarrow \xi\cdot F$ labeled by $F^{-1}$. Applying this pasting process to all boundary points on $\mathcal{T}_{\gamma}$, the new graph is the orbital graph $\Gamma_{\gamma}$.
\end{rmk}

To prove the main theorem, we will also need the family of inverse semigroups (or groups) defined by the same tile inflation process, denoted $\{G_0\}_{\mathsf{B}}$ (or $\{G\}_{\mathsf{B}}$). To define this family, we need the notion of \textit{automata}.    

\begin{defn}(\cite[Subsections 2.3.3-2.3.5]{nek22}) \label{nondetaut}
     Let $X_1,X_2\ldots $ and $X_1',X_2',\ldots $ be two sequences of \textit{alphabets} and $Q_1,Q_2\ldots $ be a sequence of \textit{sets of states}. A \textit{non-deterministic time-varying automaton} $\mathcal{A}$ consists of a sequence of transitions $T_1,T_2\ldots $ where $T_n\subset Q_n\times Q_{n+1}\times X_n\times X_n'$. The automaton is \textit{$\omega$-deterministic} if for each sequence $x_1x_2\ldots $ where $x_i\in X_i$, there exists at most one sequence of transitions of the form $(q_1,q_2,x_1,x_1'),(q_2,q_3,x_2,x_2')\ldots $ where $q_i\in Q_i$ and $x_i'=\lambda_i(q_i,x_i)$. Each initial state $q\in Q_1$ defines a transformation $X_1\times X_2\times X_3\times\ldots \rightarrow X_1'\times X_2'\times X_3'\times\ldots $. 
\end{defn}

Time-varying automata can be presented by Moore diagrams. 

\begin{defn}
Let $\mathcal{A}$ be a non-deterministic time-varying automaton. Its \textit{Moore diagram} consists of a set of vertices $Q=\bigcup_n Q_n$ and a set of edges $T=\bigcup_n T_n$, where $(q_i,q_{i+1},x,y)\in T_i$ is an arrow from $q_i\in Q_i$ to $q_{i+1}\in Q_{i+1}$ labeled by $x|y$. The state $q_{i+1}$ is an element in the \textit{section} of $q_i$ at $x$. The section (set of states) is denoted ${}_{x}{|q_i}$. 
\end{defn}

\begin{defn}\label{finconpt} Let $\mathcal{T}_{v}$ be an $n$-th level tile. An \textit{$n$-th level boundary connection} is a triple $(F,\gamma_1,\gamma_2)$ where $F\in\mathcal{S}$, and $\gamma_1,\gamma_2$ are paths of length $n$ such that $F$ is an outgoing boundary edge at $\gamma_1$ and an incoming boundary edge at $\gamma_2$. If $(F,\gamma_1,\gamma_2)$ is a boundary connection, then $(F^{-1},\gamma_2,\gamma_1)$ is also a boundary connection. 
\end{defn}

The following was proved in \cite[Proposition 5.2.19]{nek22}. 

\begin{prop} \label{Automaton} \textup{(Automaton associated with the tile inflation process)}
 Consider the following $\omega$-deterministic time-varying automaton $\mathcal{A}$. The sequence of input-output alphabets is equal to the set of the edges $E_1,E_2,\ldots $ of $\mathsf{B}$. The set of states $Q_n$ is the union of the set of trivial states $1_v$ labeled by the vertices $v\in V_n$ and the set of $n$-th level boundary connections $(F,\gamma_1,\gamma_2)$.  

For every $e\in E_n$, define a transition from the state $1_{\mathbf{s}(e)}$ to $1_{\mathbf{r}(e)}$ labeled by $e|e$.  If $e_1$ and $e_2$ are edges such that $F$ is an outgoing boundary edge at $e_1\gamma_1$ and incoming boundary edge at $e_2\gamma_2$, then for every boundary connection $(F,e_1\gamma_1,e_2\gamma_2)$, define a transition from $(F,\gamma_1,\gamma_2)$ to $(F,e_1\gamma_1,e_2\gamma_2)$ labeled by $e_1|e_2$. Otherwise, define a transition from $(F,\gamma_1,\gamma_2)$ to $1_{\mathbf{r}(e_1)}$ labeled by $e_1|e_2$. 

Then the set of initial states of the form $(F,v_1,v_2)$, for $v_1,v_2\in V_1$, defines the local homeomorphism $F$. Each non-initial non-trivial state has exactly $1$ incoming edge. 
\end{prop}

It follows that $(F,e_1\gamma_1,e_2\gamma_2)$ is an element of the section of $(F,\gamma_1,\gamma_2)$ at $e_1$. Let $F\in\mathcal{S}$ (i.e., it is defined by a set of initial states of the form $(F,v_1,v_2)\in\mathcal{A}$). Let $\gamma\in\Omega^*(\mathsf{B})$. Then the \text{section of $F$ on $\gamma$}, denoted ${}_{\gamma}{|F}$, is defined to be the set of states of the form $(F,\gamma,-)$.   

The family $\{G_0\}_{\mathsf{B}}$ is defined as follows. Each $G_i$ is generated by transformations in ${}_{\gamma}{|F}$ for all $F\in \mathcal{S}$ and $\gamma\in\Omega_i(\mathsf{B})$ on which $F$ is defined.  

\subsection{Bounded type}
\begin{defn}\label{bddtpp}
    A tile inflation process is of \textit{bounded type} if it satisfies the following conditions.
    \begin{enumerate}
   \item\label{bddtpp1} All finite tiles have uniformly bounded cardinalities of the set of boundary points and there exist only finitely many boundary points of infinite tiles. 
   \item\label{bddtpp2} The set of vertices and edges on each level of $\mathsf{B}$ have uniformly bounded cardinalities.  
    \end{enumerate}
\end{defn}

\begin{defn}\label{semibd}
    The inverse semigroup generated by $\mathcal{S}$, denoted $G_0=\langle \mathcal{S}\rangle$, is of \textit{bounded type} if it is defined by a tile inflation process of bounded type.
\end{defn}

\begin{defn}\label{grpbd}
  Let $S$ be a finite labeling set of a tile inflation process of bounded type. The inverse semigroup $G$ generated by $S$, denoted $G=\langle S\rangle$, is a \textit{group of bounded type} if,  after extending the action of some labels at boundary points of infinite tiles, the tile inflation process determines perfectly labeled orbital graphs labeled by elements in $S$. 
\end{defn}

\subsection{Incompressible elements and traverses}

\begin{defn}
     \label{trajectory} Let $v$ be a vertex in an orbital/tile graph. Let $g=F_1\ldots F_m$ for composable $F_i\in\mathcal{S}$. The \textit{trajectory} or \textit{walk} of $g$ starting at $v$ is the path \[v\cdot F_1,v\cdot F_1F_2,\ldots ,v\cdot F_1F_2\ldots F_m.\] The points $v$ and $v\cdot F_1F_2\ldots F_m$ are called \textit{initial} and \textit{final} vertices of the trajectory of $g$. 
\end{defn}

\begin{defn}\label{return} Let $G_0=\langle \mathcal{S} \rangle$ be an inverse semigroup of bounded type. Let $\xi_1,\xi_2\in\Omega(\mathsf{B})$ be two (possibly $\xi_1=\xi_2$) boundary points on the same infinite tile $\mathcal{T}_{\xi_1}$ (as a subgraph of the orbital graph $\Gamma_{\xi_1}$). Let $g=F_1\ldots F_n\in G_0$, for composable $F_1,\ldots ,F_n\in \mathcal{S}$, be an element in $G_0$. We say that $g$ is a \textit{return word} if the final vertex of the trajectory of $g$ starting at $\xi_1\cdot F_1^{-1}$ is $\xi_2\cdot F_m$, $F_1$ and $F_m$ are respectively labels of boundary edges at $\xi_1$ and $\xi_2$ on the tile $\mathcal{T}_{\xi_1}$ ($=\mathcal{T}_{\xi_2}$), while all edges in between labeled by $F_2,\ldots ,F_{m-1}$ on the trajectory of $g$ belong to $\mathcal{T}_{\xi_{1}}$ and are not boundary edges of $\mathcal{T}_{\xi_1}$. 
\end{defn}
\begin{defn}\label{incom} An element $g=F_1\ldots F_n$ is said to be \textit{incompressible} if it does not contain any return subword.
\end{defn}
\begin{defn}\label{traverse} Let $\mathcal{T}_{i,n}$ be a finite tile. Let $w=F_1F_2\ldots F_m\in \mathcal{S}^*$ be a word as above. A \textit{traverse} of $\mathcal{T}_{i,n}$ is the trajectory of a subword of $w$, denoted $w'=F_i\ldots F_j$ for $1\leq k<j \leq l$, such that the initial and final vertices of the trajectory of $w'$ are boundary points of $\mathcal{T}_{i,n}$ and all other vertices of the trajectory of $w'$ are inside $\mathcal{T}_{i,n}$ different from the boundary points.  
\end{defn}

\section{Expanding tile inflations}\label{exptile}

Expanding tile inflations were first studied for groups generated by bounded automata (acting on regular rooted trees). See \cite[Chapter V]{Bon}. In this section, we discuss more general expanding tile inflations on Bratteli diagrams.  
\begin{defn}\label{telescoping}
    Let $\mathsf{B}=((V_n)_{n=1}^{\infty}),(E_n)_{n=1}^{\infty},\mathbf{s},\mathbf{r})$ be a Bratteli diagram. Let $1=k_1<k_2<\ldots $ be an increasing sequence of integers. The \textit{telescoping} of $\mathsf{B}$ defined by the sequence is the diagram $((V'_n)_{n=1}^{\infty}),(E'_n)_{n=1}^{\infty},\mathbf{s},\mathbf{r})$ where $V'_n=V_{k_n}$, $E'_n$ is the set of paths in $\mathsf{B}$ from $V_{k_n}$ to $V_{k_{n+1}}$, and $\mathbf{s}$, $\mathbf{r}$ are the beginning and the end of the paths, respectively.  
\end{defn}

 Two Bratteli diagrams $\mathsf{B}_1$, $\mathsf{B}_2$ are \textit{equivalent}, denoted $\mathsf{B}_1\sim\mathsf{B}_2$, if and only if one can be obtained from the other by a sequence of telescopings and inverses of telescopings. If $\mathsf{B}_1\sim\mathsf{B}_2$, then $\Omega(\mathsf{B}_1)$ is naturally homeomorphic to $\Omega(\mathsf{B}_2)$. 

 \begin{defn}\label{telebounded}
     A telescoping of $\mathsf{B}$ is said to have \textit{bounded steps} if there exists $K\geq 0$ such that the sequence $\{k_n\}$ in Definition \ref{telescoping} satisfies $k_{n+1}-k_n\leq K$ for all $n\in\mathbb{N}$.
 \end{defn}
\begin{defn}
    A \textit{block} of boundary/connecting points on a finite tile $\mathcal{T}_{v_n}$ with $n\geq N$, for some $N\in\mathbb{N}$, is a subset $\mathcal{B}_n$ of the union of the set of boundary points and the set of connecting points of $\mathcal{T}_{v_n}$ in which either all continuations of elements in $\mathcal{B}_n$ are non-boundary points or for some successive embeddings 
    \begin{center}
    \begin{tikzcd}
\mathcal{T}_{v_n} \arrow[r, "\phi_{e_n}"] & \mathcal{T}_{v_{n+1}} \arrow [r, "\phi_{e_{n+1}}"] & \ldots \arrow [r, "\phi_{e_{n+m-1}}"]& \mathcal{T}_{v_{n+m-1}} \arrow [r, "\phi_{e_{n+m}}"] & \ldots,
\end{tikzcd}
\end{center}
   the sets $\mathcal{B}_{n+m}\neq \varnothing$, whose elements are obtained by appending $e_{n+m}$ to the end of some elements in $\mathcal{B}_{n+m-1}$. Here either $m\in\{1,...,M\}$ for some $M<\infty$ or $m\in\mathbb{N}$. 
\end{defn}
\begin{rmk}
In the definition above, the condition $\mathcal{B}_{n+m}\neq \varnothing$ means that there exist $\gamma_1,\gamma_2\in\mathcal{B}_{n+m-1}$ such that they are simultaneously continued by the same edge $e_{n+m}$ to be boundary/connecting points on the next level, i.e., $e_{n+m}\gamma_1,e_{n+m}\gamma_2\in\mathcal{B}_{n+m}$. It is obvious that the distances $d_{n+m-1}(\gamma_1,\gamma_2)=d_{n+m}(e_{n+m}\gamma_1,e_{n+m}\gamma_2)$, where $d_i$ are metrics for $\mathcal{T}_{v_i}$ counting the number of edges between two vertices modulo multiple edges and loops. Hence the the distances between each pair of elements in each $\mathcal{B}_k$, for $k>N$, are uniformly bounded above by the maximal distance between elements in $\mathcal{B}_N$. There is a smallest $N\in\mathbb{N}$ such that the distances between these boundary points stop grow. We call such $N$ the \textit{starting level}.
\end{rmk}
Notice that there might be a boundary point $\gamma_{k+1}$ on $\mathcal{T}_{v_{k+1}}$ not coming from continuations of elements in $\mathcal{B}_k$  whose distances to some elements in $\mathcal{B}_{k+1}$ are bounded above by by the maximal distance between elements in $\mathcal{B}_N$. In this case, we may assume $\gamma_{k+1}$ cannot be continued to be boundary points within uniformly bounded distances from elements in $\mathcal{B}_{k+m}$ for all $m> 1$. Otherwise, we may increase the starting level $N$, and include a continuation of $\gamma_{k+1}$ to be an element of $\mathcal{B}_{N}$. 

\begin{exmp}
  Let $\xi_1\neq\xi_2$ be two boundary points on the same infinite tile. Then there exists $N>0$ such that for all $n\geq N$, their $n$-th truncations belong to the same block. 
\end{exmp}
\begin{defn}
    A tile inflation process on $\mathsf{B}$ is \textit{expanding} if, after transforming $\mathsf{B}$ to $\mathsf{B}'$ by telescoping of bounded steps, for all $n\geq M$ and for each $\mathcal{T}_{i,k_{n+1}}$, every pair of blocks of boundary points of $\mathcal{T}_{i,k_{n+1}}$ come from different embeddings of tile(s) on level $k_n$.

\end{defn}

\subsection{Exponential growth of diameters}
\begin{lemma}\label{expbdy}
   Let $\mathsf{B}$ be a Bratteli diagram whose sets of edges $E_n$ satisfy $2\leq |E_n|\leq E$ for some $E>2$. Let $\{\mathcal{T}_n\}$ be any sequence of finite tiles on levels $n$ of $\mathsf{B}$ obtained by the maps $\phi_{v_n,e_n}$ in Definition \ref{inftile}. If the tile inflation process is expanding, then there exists $a>1$ such that $\max \{d(\gamma_{1,n},\gamma_{2,n})\}\sim a^n$, where $\gamma_{1,n},\gamma_{2,n}$ are boundary points on $\mathcal{T}_n$. 
\end{lemma}
\begin{proof}
   Without loss of generality, we assume that each block is a singleton and no telescoping is performed so that the argument starts from level $1$ of $\mathsf{B}$ and does not skip levels of $\mathsf{B}$. Let $\xi_{1,n},\xi_{2,n}$ be the boundary points on $\mathcal{T}_n$ representing $\max \{d(\gamma_{1,n},\gamma_{2,n})\}$. Then a geodesic $l$ connecting $\xi_{1,n},\xi_{2,n}$ is of the form $l_1F_1l_2F_2\ldots l_mF_ml_{m+1}$ where where each $l_k$ lies in some isomorphic copy of $\mathcal{T}_{j,n-1}$, for various $j\in\{1,2,\ldots ,|V_n|\}$, and $F_k\in \mathcal{S}$ is the label of a boundary edge of $\mathcal{T}_{j,n-1}$. Let $p_k,q_k$ be respectively the beginning and ending vertices of $l_k$. Note that $p_1=\xi_{1,n}$ and $q_{m+1}=\xi_{2,n}$ and other $p_k,q_k$ are connecting points of $\mathcal{T}_{n+1}$. It follows that 
   \[d(\xi_{1,n},\xi_{2,n})=d(\xi_{1,n},p_1)+\sum\limits_{k=2}^{m}d(p_k,q_k)+d(p_{m+1},\xi_{2,n})+m.\]
   Hence \[d(\xi_{1,n},\xi_{2,n})\geq d(\xi_{1,n},p_1)+d(p_{m+1},\xi_{2,n}).\addtag\label{ieq1}\] Relabel $p_{1,n}=p_1$ and $p_{2,n}=p_{m+1}$. The inequality above can be rewritten as
   \[d(\xi_{1,n},\xi_{2,n})\geq d(\xi_{1,n},p_{1,n})+d(p_{2,n},\xi_{2,n}).\addtag\label{ieq2}\]
   Since the inflation process is expanding, we can write $e_1\xi_{1,n}=\xi_{1,n-1}$, $e_{1,n}p_{1,n}=p_{1,n-1}$, $p_{2,n}=e_{2,n}p_{2,n-1}$, and $\xi_{2,n}=e_{2,n}\xi_{2,n-1}$, for $e_{1,n}\neq e_{2,n}$. Then the vertices $\xi_{1,n-1}$, $p_{1,n-1}$, $p_{2,n-1}$ and $\xi_{2,n-1}$ are boundary points of tiles on level $n-1$, where $\xi_{1,n-1},p_{1,n-1}$ belong to the same tile and $p_{2,n-1},\xi_{2,n-1}$ belong to the same tile. Since the inflation process is expanding, $\xi_{1,n-1},p_{1,n-1}$ come from different embeddings of tile(s) on level $n-2$, and $p_{2,n-1},\xi_{2,n-1}$ come from different embeddings of tile(s) on level $n-2$. By the same reasoning we have 
   \[d(\xi_{1,n-1},p_{1,n-1})\geq d(\xi_{1,n-1},p_{3,n-1})+d(p_{4,n-1},p_{1,n-1}) \addtag\] and
   \[d(p_{2,n-1},\xi_{2,n-1})\geq  d(p_{2,n-1},p_{5,n-1})+d(p_{6,n-1},\xi_{2,n-1}).\addtag\]
The above inequalities hold for all $n>1$. It follows that $d(\xi_{1,n},\xi_{2,n})\geq 2^n$. 

On the other hand, since each $|E_n|\leq E$, $d(\xi_{1,n},\xi_{2,n})\leq E^n$. It follows that
\[2^n\leq\max\{d(\xi_{1,n},\xi_{2,n})\}\leq E^n.\addtag \label{exppbounds}\]

To determine the exact exponential growth rate, let $d_n=\max\{d(\xi_{1,n},\xi_{2,n})\}$ and $s_n=\log d_n$. We \textbf{claim} that the sequence $\{d_n\}$ is \textit{sub-multiplicative}, i.e., for all $m,n\in\mathbb{N}$, there exists a constant $C\geq 1$ such that 
\[d_{m+n}\leq C\cdot d_md_n.\addtag \label{submul}\]
This follows from (\ref{exppbounds}) and the assumption that the set of edges of $\mathsf{B}$ has uniformly bounded cardinalities. Indeed, for each $m$, $2^m\leq d_m\leq E^m$ and thus $2^{m+n}\leq d_md_n\leq E^{m+n}$. Since the tile inflation process is of bounded type, there exists constants $1<c_1<c_2$ such that $d_{n+1}\leq c_1\cdot d_n\leq c_22^n$. Hence $d_{m+n}\leq c_22^{m+n}$. It follows that 

\[\dfrac{d_{m+n}}{d_md_n}\leq \dfrac{c_22^{m+n}}{2^m2^n}=c_2.\addtag\label{ccc2}\]

Let $C=c_2$ and the claim is proved.

It then follows from (\ref{submul}) that the sequence $\{s_n\}$ is \textit{sub-additive}, i.e., 
\[s_{m+n}\leq s_m+s_n +\log C. \addtag\label{subadd} \]
By Fekete's Subadditive Lemma (see, for instance, \cite{fek23}) and the fact that $\lim\limits_{n\rightarrow\infty}\frac{\log C}{n}=0$, we have that the limit \[\lim\limits_{n\rightarrow \infty}\dfrac{s_n}{n}\]
exists and is denoted $\alpha>0$. Hence $d_n$ is asymptotically equivalent to $\exp(\alpha n)$. Since $2\leq d_n\leq E$, we have $\log 2\leq\alpha\leq \log E$. Let $a=\exp(\alpha)$ and the proof is complete. 
\end{proof}

The following was proved in \cite[Proposition 3.9]{kua24a}.
\begin{prop}\label{expdiam}
    Let $\{\mathcal{T}_{n}\}$ a the sequence of tiles indexed by the levels $n$, for $n\in\mathbb{N}$. Then \textup{Diam}$(\mathcal{T}_{i,n})$ is asymptotically equivalent to $d(\gamma_{1,n},\gamma_{2,n})$ where $\gamma_{1,n},\gamma_{2,n}$ are boundary points on $\mathcal{T}_{i,n}$ with $\{d(\gamma_{1,n},\gamma_{2,n})\}$ growing exponentially. Consequently, the sequence of diameters of a tile inflation has exponential growth.    
\end{prop}
By Lemma \ref{expbdy} and Proposition \ref{expdiam}, we have the following.

\begin{cor}\label{equia}
    Let $\mathsf{B}$ be a Bratteli diagram whose sets of edges $E_n$ satisfy $2\leq |E_n|\leq E$ for some $E>2$. Then each tile on level $n$ of $\mathsf{B}$ obtained by an expanding tile inflation process has diameter asymptotically equivalent to an exponential function $a^n$, for some $a>1$. 
\end{cor}

\subsection{Linearly repetitive actions}
\begin{defn}
Let $\mathcal{T}$ be a finite tile and $\Gamma$ be an orbital graph or infinite tile. The tile $\mathcal{T}$ is said to be \textit{linearly repetitive} on $\Gamma$ if there exists a constant $K>1$ and given any $v\in V(\Gamma)$, there exists an isomorphic copy $\mathcal{T}'$ of $\mathcal{T}$ that is contained in $B_v(Kr)$, where $B_v$ is the ball on $\Gamma$ and $r$ is the diameter of $\mathcal{T}$. The tile inflation process and the action $\Omega(\mathsf{B})\curvearrowleft G_0$ are called \textit{linearly repetitive} if all finite tiles are linearly repetitive on all infinite tiles/orbital graphs. 
\end{defn}
\begin{prop}\label{linrep}
    Let $G_0$ be an inverse semigroup satisfying Condition \ref{bddtpp1} in Definition \ref{bddtpp} and acting minimally on $\Omega(\mathsf{B})$. Let $\mathcal{T}_1,\mathcal{T}_2,\ldots ,\mathcal{T}_n\ldots $ be any sequence of tiles where $\mathcal{T}_n$ is a tile on level $n$ of $\mathsf{B}$. Suppose their diameters grow exponentially, i.e., $\textup{Diam}(\mathcal{T}_n)\sim a^n$ for some $a>1$. The following are equivalent. 
    \begin{enumerate}[label=(\arabic*)]
        \item Finite tiles are linearly repetitive on infinite tiles.\label{linone}
        \item The set of vertices and edges on each level of $\mathsf{B}$ have uniformly bounded cardinalities, and for each $n$, the smallest $l\in\mathbb{N}$ (depending on $n$) such that all tiles on level $n+l$ contain all tiles on level $n$ is uniformly bounded. \label{lintwo}
    \end{enumerate}
\end{prop}
Notice that since $\mathsf{B}$ is simple, the assumption that the sets of edges of $\mathsf{B}$ have uniformly bounded cardinality implies the assumption that the sets of vertices of $\mathsf{B}$ have uniformly bounded cardinality.  
\begin{proof}
   (\ref{lintwo}$\implies$\ref{linone}):  Suppose, for each level of $\mathsf{B}$, the numbers of edges and vertices are bounded by $M>2$ and $N>1$, respectively, and $l$ is bounded by $L$. Let $n\in\mathbb{N}$ be any number. By minimality, there exists the smallest $l(n)>0$ such that all tiles on level $n$ of $\mathsf{B}$ are contained as isomorphic copies in each tile on level $m=n+l(n)$. For the same reason, all tiles on level $m$ are contained as isomorphic copies in each tile on level $m+l(m)$ for the smallest $l(m)$.  Fix $\mathcal{T}_{j,m+l(m)}$. It contains at most $M^{l(m)}$ copies of tiles on level $m$ and at most $M^{l(m)+l(n)}\leq M^{2L}$ copies of tiles on level $n$. Since the diameters of all tiles on level $n$ are asymptotically equivalent to $a^n$ for some $a>1$, fixing $i$, the gap between consecutive appearance of $\mathcal{T}_{i,n}$ in $\mathcal{T}_{j,m+l(m)}$ has diameter not more than $(2M^L-2)a^n\sim (2M^L-2)\cdot\textup{Diam}(\mathcal{T}_{i,n})$. This is because each $\mathcal{T}_{j,m}$ contains at least one copy of $\mathcal{T}_{j,n}$, and if $\mathcal{T}_{j_1,m}$, $\mathcal{T}_{j_1,m}$ are adjacent to each other on $\mathcal{T}_{j,m+l(m)}$, then the gap between $\mathcal{T}_{i,n}$ in $\mathcal{T}_{j_1,m}$ and $\mathcal{T}_{i,n}$ in $\mathcal{T}_{j_2,m}$ has diameter not more than the sum of the diameters of all tiles on level $n$ contained in them. Since infinite tiles are obtained from finite tile inflations, the statement is true for any $\mathcal{T}_{i,n}$ in any infinite tiles. 

   (\ref{linone}$\implies$\ref{lintwo}): Let $M(n)$, $N(n)$ be the total numbers of edges and vertices on each level of $\mathsf{B}$, respectively, and let $l(n)$ be as above. Suppose first that $N(n)$ is not uniformly bounded. This implies that $M(n)$ is not uniformly bounded. Suppose, without loss of generality, that $N(n)$ is monotonously increasing. For each $n$, suppose $l(n)$ is bounded by $L$. We may assume each $l(n)=1$, which means that each tile on level $n+1$ is covered by all tiles on level $n$, for all $n\in\mathbb{N}$. Let $\mathcal{T}_{i,n}$ be a tile on level $n$. Then the largest gap between consecutive appearances of $\mathcal{T}_{i,n}$ on level $n+2$ will have diameter bounded below by $(N(n)-k(n))\cdot \textup{Diam}(\mathcal{T}_{i,n})$, where $k(n)$ is a bounded number. Since the function $N(n)$ is increasing, we conclude that finite tiles cannot be linearly repetitive. 
   
   Now suppose $l(n)$ is not bounded but $N(n)$ is bounded. We may also assume without loss of generality that $l(n)$ is monotonously increasing and $N(n)\geq 2$ is a constant. Then for each $\mathcal{T}_{i,n}$, the largest gap between consecutive appearances of it on a tile on level $n+l(n)+l(n+l(n))$ is bounded below by $E^{l(n)}\cdot \textup{Diam}(\mathcal{T}_{i,n})$ where $E$ is a constant greater than $1$. Since $l(n)$ is increasing, some finite tiles are not linearly repetitive. 
   \end{proof}
\begin{rmk}
By the proposition above, linear repetitivity of finite tiles is a stronger condition than bounded type. However, even if $l(n)$ above is not uniformly bounded, for each level $n$ of $\mathsf{B}$, we can perform a telescoping $\mathsf{B}'$ of $\mathsf{B}$ such that all tiles on level $n+l'(n)$ contain all tiles on level $n$ for uniformly bounded $l'$ defined on $\mathsf{B}'$. Hence we may always assume that $l(n)$ is uniformly bounded to ensure linear repetitivity. Informally, linear repetitivity is equivalent to bounded type up to telescopings. Therefore, for the rest of this paper, we will assume the action of $G_0$ on $\Omega(\mathsf{B})$ is linearly repetitive. 
\end{rmk}

\subsection{Polynomial growth of orbital graphs} 

For a finite tile $\mathcal{T}_n$, denote by $|\mathcal{T}_n|:=\#V(\mathcal{T}_n)$ its \textit{cardinality} of vertices.
\begin{lemma}\label{equidiacar}
  The cardinalities of finite tiles on level $n$ of $\mathsf{B}$ with $2\leq E_n\leq E$
  obtained by an expanding tile inflation process of bounded type are asymptotically equivalent to $b^n$ for some $b>1$. 
 \end{lemma}

 \begin{proof}
It is enough to show that the sequence $\{|\mathcal{T}_n|\}_{n\in\mathbb{N}}$ for any tile $\mathcal{T}_n$ on level $n$ is sub-multiplicative. Indeed, we only need to replace the $d_m,d_n$ in (\ref{exppbounds}), (\ref{submul}), (\ref{ccc2}) with $|\mathcal{T}_m|,|\mathcal{T}_n|$. The rest of the argument is the same as in the proof of Lemma \ref{expbdy}.
 \end{proof}

\begin{prop}
     Let $\mathsf{B}$ be a Bratteli diagram whose sets of edges $E_n$ satisfy $2\leq |E_n|\leq E$ for some $E>2$. Let $G_0$ be defined by an expanding and linearly repetitive tile inflation process on $\mathsf{B}$. Let $\gamma\in\Omega(\mathsf{B})$ be any point and let $\mathcal{T}_{\gamma}$ be the infinite tile containing $\gamma$. Let $B_{\gamma}(R)\subset\mathcal{T}_{\gamma}$ be a ball of radius $R$ centered at $\gamma$. Then there exist $C>1$ and $C'>0$ such that the cardinality $\# B_{\gamma}(R)$ satisfies 
     \[C'R^d\leq\# B_{\gamma}(R)\leq CR^d.\addtag\label{polyball}\] 
\end{prop}
\begin{proof}
    The ball $B_{\gamma}(R)$ is covered by a uniformly bounded number $M$ of isomorphic copies of finite tiles on level $n$ whose diameters are asymptotic equivalent to $a^n$ for some $a>1$ (see Corollary \ref{equia}). In order for $a^n$ to be asymptotic equivalent to $R$, denoted $a^n\sim R$, we must have $n\sim \frac{\log R}{\log a}$. Write $n=k\log_{a}R+c$ for $k>0$. Hence \[\#B_{\gamma}(R)\leq M|\mathcal{T}_n|^{k\log_aR+c}\sim Mb^{k\log_aR+c},\] by Lemma \ref{equidiacar}, while $b^{k\log_aR+c}$ is asymptotically equivalent to $R^d$ for some $d>0$. Since $B_{\gamma}(R)$ is covered by $M$ tiles with cardinality bounded above by $M'R^d$ for some $M'>0$, taking $C=M'M$, the upper bound is proved. 
    
    Now let us show the lower bound. Since the tile inflation is linearly repetitive, for any finite tile $\mathcal{T}_n$, there exists $K>1$ such that the ball $B_{\gamma}(K\cdot\textup{Diam}(\mathcal{T}_n))$ contains an isomorphic copy of $\mathcal{T}_n$ with $\textup{Diam}(\mathcal{T}_n)\sim a^n$. Let $R\sim a^n$, and write $n=k\log_aR+c$ as above. Then 
    \[\#B_{\gamma}(KR)\geq M''b^n=M''b^{k\log_aR+c}\sim M''R^d\] 
    for the same $d$ as above and for some $M''>0$. This implies that 
    \[\#B_{\gamma}(R)\geq \frac{1}{K}c'M''R^d\] 
    for some $c'>0$. Taking $C'= \frac{1}{K}c'M''$, the lower bound is proved.
\end{proof}
\begin{rmk}
Linear repetitivity is only used to prove the lower bound. 
\end{rmk}

 Since, for a tile inflation process of bounded type, the orbital graphs and infinite tile graphs differ for at most a bounded cardinality of boundary edges, we have the following. 
\begin{cor}
    The orbital graphs of an inverse semigroup of bounded type have polynomial growth given by (\ref{polyball}). 
\end{cor}

\section{Proof of the theorem}\label{prvmain}
\begin{defn}\label{orbequii}
    Let $G_1,G_2$ be inverse semigroups acting on $\Omega(\mathsf{B})$. The action of $G_1$ is \textit{orbit equivalent} to the action of $G_2$ if there exist a pair of homeomorphisms $\varphi,\psi:\Omega(\mathsf{B})\rightarrow \Omega(\mathsf{B})$ such that for any $\gamma\in \Omega(\mathsf{B})$ and any $g_1\in G_1$ defined at $\gamma$, there exists $g_2\in G_2$ defined at $\varphi(\gamma)$ such that $\varphi(\gamma\cdot g_1)=\varphi(\gamma)\cdot g_2$, and for any $g_2'\in G_2$ defined at $\gamma$, there exists $g_1'\in G_1$ defined at $\psi(\gamma)$ such that $\psi(\gamma\cdot g_2')=\psi(\gamma)\cdot g_1'$. 
\end{defn}

\begin{defn}\label{strobequi}
    The actions of $G_1,G_2$ on $\Omega(\mathsf{B})$ above are said to be \textit{strongly orbit equivalent} if the relations $\varphi(\gamma\cdot g_1)=\varphi(\gamma)\cdot g_2$ and $\psi(\gamma\cdot g_2')=\psi(\gamma)\cdot g_1'$ in Definition \ref{orbequii} hold in a neighborhood of $\gamma$ for any $\gamma\in\Omega(\mathsf{B})$. 
\end{defn}
 Recall that the \textit{topological full group} of $G_0$, denoted $\mathsf{F}(G_0)$, is the set of all automorphisms $f:\Omega(\mathsf{B})\rightarrow \Omega(\mathsf{B})$ defined by the rule $\gamma\cdot f=\gamma\cdot g$ for each $\gamma\in \Omega(\mathsf{B})$ and for some $g\in G_0$ defined in a neighborhood of $\gamma$, i.e., every element of $\mathsf{F}(G_0)$ acts ``piecewise" as an element of $G$. The following is the main result. 
 \begin{thm}\label{main}
     Let $\mathsf{B}$ be a simple Bratteli diagram. Suppose the inverse semigroup $G_0=\langle \mathcal{S}\rangle$ is defined by an expanding and linearly repetitive tile inflation process and acts minimally on $\Omega(\mathsf{B})$. Let $G<\mathsf{F}(G_0)$ be a finitely generated subgroup that is orbit equivalent to $G_0$. If the set of incompressible elements of $G_0$ is finite, then there exists $\alpha\in (0,1)$ such that the growth function $\gamma_G$ of $G$ satisfies $\gamma_G(R)\preccurlyeq \exp(R^{\alpha})$. 
 \end{thm}
The proof will follow from several propositions.  

\begin{prop}\label{orbequi}
   If $G<\mathsf{F}(G_0)$ is orbit equivalent to $G_0$, then it is strongly orbit equivalent to $G_0$. \end{prop}
\begin{proof}
    This follows directly from the definition of topological full groups. Indeed, for any $\gamma\in \Omega(\mathsf{B})$, any $g\in \mathsf{F}(G_0)$ acts as an element in $G_0$ in a neighborhood of $\gamma$. 
\end{proof}

Recall two infinite paths $\gamma_1,\gamma_2\in\Omega(\mathsf{B})$ are called \textit{cofinal} if they differ for finitely many edges.  

 \begin{cor}
    If $G<\mathsf{F}(G_0)$ is finitely generated and orbit equivalent to $G_0$, then $G$ is a group of bounded type. 
\end{cor}
\begin{proof}
   By Proposition \ref{orbequi}, $G$ is strongly orbit equivalent to $G_0$. Hence the boundary points on infinite tiles of $G$ are cofinal to those of $G_0$. Since $G$ is finitely generated, it follows that the number of boundary points on infinite tiles of $G$ is finite. Now let $\mathcal{T}$ be a finite tile of $G_0$ and $B_{\mathcal{T}}$ be its set of boundary points. Then there exists $B>0$ such that $|B_{\mathcal{T}}|\leq B$. Hence the number of boundary points of any finite tile of $G$ is bounded above by $B+|S|$, where $S$ is the (finite) generating set of $G$. Hence $G$ is of bounded type.
 \end{proof}

The following was proved in \cite[Proposition 4.3]{kua24a}, which is a modification of \cite[Proposition 2.7]{BNZ}.

\begin{prop} \label{ultesti}
 Suppose that there are $t_0>1$ and a function $\Psi(t)$ such that $F_{w}(t)\leq |w|\cdot\Psi(t)$ for all $w\in S^*$, all $t<t_0$. Set $\beta=\limsup_{n\rightarrow \infty}(\sum_i|\mathcal{T}_{i,n}|)^{1/n}$. Then the growth function $\gamma_G(R)$ of $G$ satisfies $$\gamma_G(R)\preccurlyeq \exp(R^{\alpha})$$ for all $\alpha>\frac{\log \beta}{\log\beta+\log t_0}$. 
\end{prop}

We also need the following technical definition. 
\begin{defn}\label{reduced}
    Let $w=s_1\ldots s_m\in S^*$ be a word. We say that $w$ is \textit{reduced} if all $s_i\neq Id$ and either all consecutive $s_i,s_{i+1}$ are labels of boundary edges at different boundary points on infinite tiles or at least one of them is finitary. Otherwise, the word $w$ is said to contain \textit{non-reduced subwords}.  
\end{defn}

\begin{proof}[Proof of Theorem \ref{main}] Let $n$ be any level of $\mathsf{B}$. Then all the tiles on level $n+l(n)$ contain (isomorphic copies of) all the tiles on level $n$, for uniformly bounded $l(n)$. Hence, by applying a telescoping of bounded steps of $\mathsf{B}$, we may assume, without loss of generality, that all the tiles on level $n+1$ contain all the tiles on level $n$ and blocks of boundary points on level $n+1$ come from different embeddings of tile(s) on level $n$. 
Let $\mathcal{T}_{i,n}$ be a tile on level $n$. Then one of the following cases must happen. 
\begin{enumerate}
    \item \label{stranded} For any embedding of $\mathcal{T}_{i,n}$ into level $n+1$, any traverse of isomorphic copies of $\mathcal{T}_{i,n}$ induced by traverses of tiles on level $n+1$ contains a return. In this case $\mathcal{T}_{i,n}$ is called \textit{stranded}. 
  \item \label{blockinherit} There exists an embedding of $\mathcal{T}_{i,n}$ into $\mathcal{T}_{j,n+1}$ such that one of the blocks of boundary points of $\mathcal{T}_{i,n}$ becomes a block of boundary points of $\mathcal{T}_{j,n+1}$.
  \item \label{inmiddle} For all embeddings of $\mathcal{T}_{i,n}$, at least $2$ blocks of boundary points of the isomorphic copies of $\mathcal{T}_{i,n}$ are connected to other embeddings of tile(s) while no isomorphic copy of $\mathcal{T}_{i,n}$ contains any block of boundary points of level $n+1$. 
\end{enumerate}
 
The rest of the proof is similar to that of \cite[Theorem 4.1]{kua24a}. Fix $w=s_1\ldots s_m\in S^*$. Let $\mathcal{T}_{j,n+1}$ be a tile on level $n+1$. Let $\tau$ be a traverse of $\mathcal{T}_{j,n+1}$ (viewed as a subword of $w$). Then $\tau$ induces a finite sequence of traverses of tiles on the $n$-th level. Denote this sequence by $\Phi(\tau)=(\tau_1,\tau_2,\ldots ,\tau_k)$.  Let $\mathbf{w}=x_1x_2\ldots x_k$ be the sequence of labels of the incoming edges of each induced traverse at each boundary point. Each $x_i$ is a product of elements in $S$, and $x_1$ can be any edge going into the tile $\mathcal{T}_{j,n+1}$. Let $\Theta_n$ be the set of all traverses of tiles on level $n$, for $n\in\mathbb{N}$. We may assume, without loss of generality, that there are no stranded tiles. 

We distinguish $2$ cases of $\mathbf{w}$: 1. $\mathbf{w}$ is a reduced word; 2. $\mathbf{w}$ contains non-reduced subwords. 

Suppose first the word $\mathbf{w}$ is reduced. 
\begin{defn}\label{lastmoment}
    Define a map $\phi:\Theta_{n+1}\rightarrow\Theta_{n}$ by letting $\phi(\tau)$ be the \textbf{last} traverse induced by $\tau$ of an isomorphic copy of $\mathcal{T}_{i,n}$ that satisfies Case \ref{blockinherit} above, where $\tau$ is a traverse starting at the block of boundary point of $\mathcal{T}_{i,n}$ coinciding with that of $\mathcal{T}_{j,n+1}$. Call the map $\phi$ the \textit{last-moment map}, and $\phi(\tau)$ the \textit{last induced traverse}.
\end{defn}
Note that the traverse $\tau$ of $\mathcal{T}_{j,n+1}$ is mapped to a union of sets of traverses of different tiles on level $n$. Whenever the sequence $\mathbf{w}=x_1\ldots x_k$ contains a return subword, denoted $\mathbf{w}'=x_{l_1}\ldots x_{l_2}$, the trajectory of $\mathbf{w}'$ must traverse isomorphic copies of the same tile. denoted $\mathcal{T}_{l,n}$, on level $n$ at the beginning $x_{l_1}$ and the end $x_{l_2}$. Hence this trajectory can be moved to the infinite tile that defines this return, which will traverse the same isomorphic copy of $\mathcal{T}_{l,n}$ at the beginning and the end. It follows that there exists an embedding of $\mathcal{T}_{l,n}$ to level $n+1$ that satisfies Case \ref{blockinherit} above. Hence the set of traverses of $\mathcal{T}_{l.n}$ has a nontrivial intersection with the range of $\phi$.  

%For $j_1\neq j_2$, we require that the traverses of $\mathcal{T}_{j_1,n+1}$ and $\mathcal{T}_{j_2,n+1}$ are mapped to unions of different sets of traverses. 

The map $\phi$ is injective since the trajectory of the induced traverse on $\mathcal{T}_{i,n}$ uniquely determines the traverse of $\mathcal{T}_{j,n+1}$ that induces it. We would like to show $\phi$ is non-surjective. 

Now let $\tau$ be a traverse of any tile on the $(n+1)$-st level. Write \[\Phi(\tau)=(\tau_1,\ldots ,\tau_j,\tau_{j+1},\ldots ,\tau_k)\] where $\tau_j=\phi(\tau)$. Consider the sub-label $\mathbf{w}'=x_jx_{j+1}\ldots x_l$ of $\mathbf{w}$, for $j< l \leq k$. If $\mathbf{w}'$ itself is a return word, then $\tau_j$ cannot be in the image of $\phi$ since $\tau_{l}$ is the last one. Now suppose $\mathbf{w}$ is long enough (with length greater than $M$). Then by the assumption, $\mathbf{w}$ contains a return subword. Let $x_{l_1}\ldots x_{l_2}$ represent this subword, and denote the corresponding sequence of traverse by $(\tau_{l_1},\ldots ,\tau_{l_2})$. Note that the trajectory of the induced traverses might not return to the same vertex, but if we move this subword to the corresponding tile that defines this return, its trajectory will return to the same vertex on that tile. Then the traverse $\tau_{l_1}$ cannot be in the image of $\phi$ for the same reason. It follows that for each $M$ induced traverses, there is at least $1$ that is not in the image of $\phi$. This implies that $$T_{n+1}\leq \dfrac{M-1}{M}T_n.$$ Then the generating function $F_w(t)=\sum\limits_{n=1}^{\infty}T_nt^n$ satisfies
\[F_w(t)-T_0\leq \frac{M-1}{M}tF_w(t),\] \[\implies F_w(t)\leq\dfrac{T_0}{1-\frac{M-1}{M}t}.\addtag\label{redd}\] When $t< \frac{M}{M-1}$, the right hand side is positive. Note that $T_0=|w|$. Hence we can take $\Psi(t)=\frac{1}{1-\frac{M-1}{M}t}$.

%By Proposition \ref{ultesti}, the growth function $$\gamma_G(R)\preccurlyeq \exp(R^{\alpha})$$ for every $\alpha>\dfrac{\log(\beta)}{\log(\beta)+\log(\frac{M}{M-1})}$, where $\beta=\lim\sup\limits_{n\rightarrow \infty}(\sum_i|\mathcal{T}_{i,n}|)^{1/n}$. 

Now suppose $\mathbf{w}$ contains non-reduced subwords. Let $n\in\mathbb{N}$ be big enough. Let $r<n$ be an integer satisfying the following properties. For any traverse $\tau$ of $\mathcal{T}_{i,n+1}$ as above, denote by $\Phi(\tau)=(\tau_1,\ldots ,\tau_k)$ the sequence of traverse of tiles on level $r$ of $\mathsf{B}$. Let $\mathbf{w}_r=x_1\ldots x_k$ be the sequence of labels of incoming edges at each boundary point, where each $x_i$ is the section of some $s\in S$ on a boundary point or connecting point $\gamma\in\Omega_{r-1}(\mathsf{B})$, and $\mathbf{w}_r$ is a reduced word. In other words, the word $\mathbf{w}_r$ is an element in the group $G_r\in\{G\}_{\mathsf{B}}$. We can choose $r$ such that no consecutive $x_i,x_{i+1}$ are labels of boundary edges (of $G_r$). For example, if $x_i$ is the label of a boundary edge, then $x_{i+1}$ is the label of an edge at a connecting point, which is finitary. Define the last moment map $\phi:\Theta_{i,n+1}\rightarrow\Theta_{i,r}$ in the same way as Definition \ref{lastmoment}. By the assumption, every subword of $\mathbf{w}_r$ of length $M$. By a similar argument as above, for each $M$ induced traverses, there is at least $1$ that is not in the image of $\phi$. This implies that 
\[T_{n+1}\leq\dfrac{M-1}{M}T_r.\]
The generating function $F_w(t)$ satisfies
\[Mt^{-(n+1-r)}(F_w(t)-(T_0+T_1t+\ldots +T_{n-r}t^{n-r}))\leq (M-1)F_w(t)\] 
\[\implies F_w(t)\leq |w|\dfrac{g(t)}{t^{-(n+1-r)}-\frac{M-1}{M}},\addtag\label{nonredd}\]
where $g(t)$ is positive. Notice that $n-r$ is uniformly bounded. Write $n_1\leq n-r\leq n_2$ for some $n_1,n_2\in\mathbb{N}$.  Taking $\Psi(t)=\dfrac{g(t)}{t^{-(n+1-r)}-\frac{M-1}{M}}$, then $F_w(t)\leq|w|\cdot\Psi(t)$ whenever $t<\sqrt[n+1-r]{\frac{M}{M-1}}$. Take 
\[t<\sqrt[n_2+1]{\frac{M}{M-1}}\leq \sqrt[n+1-r]{\frac{M}{M-1}}.\addtag\label{nonreddt}\] 

Combining (\ref{redd}), (\ref{nonredd}) and (\ref{nonreddt}), and by Proposition \ref{ultesti}, we conclude that the growth function
 \[\gamma_G(R)\preccurlyeq \exp(R^{\alpha})\] for every
 \[\alpha>\max\left\{\dfrac{\log(\beta)}{\log(\beta)+\log(\frac{M}{M-1})},\dfrac{\log(\beta)}{\log(\beta)+\log\left(\sqrt[n_2+1]{\frac{M}{M-1}}\right)}\right\},\]
 where $\beta=\lim\sup\limits_{n\rightarrow \infty}(\sum_i|\mathcal{T}_{i,n}|)^{1/n}$. The proof is complete.
\end{proof}

  \begin{prop}\label{torsion}
       Let $G$ be a finitely generated group with finitely many incompressible elements acting minimally on $\Omega(\mathsf{B})$. Then $G$ is periodic.   
  \end{prop} 
  \begin{proof}
The idea of the proof is similar to that of \cite[Theorem 4.1]{nekrash18}. Since the set of incompressible elements of $G$ is finite, there exists $M>0$ such that any $g\in G$ with $|g|\geq M$ contains a return. Fix a non-identity $g\in G$. Let us first prove the following lemma.
\begin{lemma}
Fix any finite tile $\mathcal{T}_g$ whose diameter is greater than or equal to $|g|$. Then for each vertex $v\in\mathcal{T}_g$ and every embedding $\phi$ of $\mathcal{T}_g$ into an infinite orbital graph there exists $k>0$ such that $\phi(v)\cdot g^k\in \phi(\mathcal{T}_g)$. 
\end{lemma}
\begin{proof}
Suppose $\mathcal{T}_g$ is on level $n$ of $\mathsf{B}$ for some $n\in\mathbb{N}$. Let $\gamma\in\Omega(\mathsf{B})$ be a point having $v$ as its prefix. Consider the trajectory \[\{\gamma\cdot g,\gamma\cdot g^2,\ldots ,\gamma\cdot g^k\} \] for some $k$ big enough. If the trajectory stays inside the isomorphic copy of $\mathcal{T}_g$ containing $\gamma$, for all $k$, then we can choose $\phi$ to be any embedding. If the trajectory returns to $\mathcal{T}_g$, then we can also choose $\phi$ to be any embedding. 

Suppose the conclusion is not true for some $\mathcal{T}_g$, i.e., for some orbital graph $\Gamma$ there exists an embedding $\phi:\mathcal{T}_g\rightarrow\Gamma$ such that the trajectory \[\{\phi(v)\cdot g,\phi(v)\cdot g^2,\ldots ,\phi(v)\cdot g^k\} \] does not come back to $\phi(\mathcal{T}_g)$, for all $k\in\mathbb{N}$. Then the length of this trajectory must approach infinity as $k\rightarrow \infty$ since otherwise $g$ will not be invertible. Hence this trajectory induces traverses of (isomorphic copies of) tiles on level $n$ of $\mathsf{B}$. Let $\mathbf{w} =x_1x_2\ldots x_j$ be
the sequence of labels of the incoming edges of each induced traverse at each boundary
point. Note that the word $\mathbf{w}$ is an element of $G$. By assumption, each subword of $\mathbf{w}$ with length at least $M$ contains a return subword. We may assume, without loss of generality, that each subword of length $M$ contains a single return subword.  Write $\mathbf{w}=x_1\ldots x_Mx_{M+1}\ldots x_{2M}\ldots x_j$. We will also use subwords of $\mathbf{w}$ to represent trajectories. Then there is an embedding $\phi_1$ of $\mathcal{T}_g$ such that the trajectory represented by $x_1\ldots x_M$ returns to a tile. Let $\mathcal{T}^{(1)}$ be an (isomorphic of) finite tile containing $\phi_1(\mathcal{T}_g)$ and the trajectory $x_1\ldots x_M$. Then there exists an embedding $\phi_2$ of $\mathcal{T}^{(1)}$ such that the trajectory $x_{M+1}\ldots x_{2M}$ returns to a tile. Let $\mathcal{T}^{(2)}$ be an (isomorphic of) finite tile containing $\phi_2(\mathcal{T}^{(1)})$ and the trajectory $x_1\ldots x_Mx_{M+1}\ldots x_{2M}$. Continue this process. Since there are only finitely many non-isomorphic returns of length not more than $M$, and $g^k$ is periodic, it follows that after successive embeddings of $\mathcal{T}^{(i)}$, the trajectory $\mathbf{w}$ contains an eventually periodic sequence of returns, for $k$ big enough. Suppose, without loss of generality, that the sequence of returns is periodic. Let $x_{l_1}\ldots x_{l_2}$, for $1\leq l_1<l_2\leq M$ be the return contained in $x_1\ldots x_M$. Let $x_{l_3}\ldots x_{l_4}$, for $iM+1\leq l_3<l_4\leq 2iM$ be the first isomorphic copy of $x_{l_1}\ldots x_{l_2}$, contained in $x_{iM+1}\ldots x_{2iM}$. Then the word $x_{l_2}\ldots x_{l_3}$ is a return word. It follows that there is an embedding $\psi$ of $\mathcal{T}^{(i-1)}$ into an orbital graph $\Gamma$ such that the trajectory $x_{l_1}\ldots x_M\ldots x_{iM+1}\ldots x_{l_3}$ is a cycle. However, this is impossible since if we apply $g^{-k}$ on $x_{l_1}\ldots x_M\ldots x_{iM+1}\ldots x_{l_4}$ on $\Gamma$, the trajectory will not touch $\psi(\mathcal{T}_g)$, contradicting the fact that $g$ is invertible. 
\end{proof}
Let $\Gamma$ be an orbital graph. Then $\Gamma$ is a disjoint union of (isomorphic copies) of all tiles on level $n$ connected at their boundary points. By the lemma above, the iteration of $g$ on any vertex $\gamma\in\Gamma$ has bounded orbits. Since $\Gamma$ is linearly repetitive, all $g$-orbits all vertices of $\Gamma$ are uniformly bounded. Hence there exists $m$ such that $g^m$ acts trivially on the vertices of $\Gamma$. Since the set of vertices of $\Gamma$ is dense in $\Omega(\mathsf{B})$, it follows that $g^n=Id$. This completes the proof.  
 \end{proof}

 By Theorem \ref{main}, Proposition\ref{torsion} and Gromov's theorem on polynomial growth \cite{Gro81}, we have the following. 

 \begin{cor}
  Let $G_0$ be an inverse semigroup of bounded type with finitely many incompressible elements. Then any finitely generated subgroup $G<\mathsf{F}(G_0)$ that is orbit equivalent to $G_0$ has intermediate growth. 
 \end{cor}
\section{Dihedral groups $D_{\infty}$}\label{dihegrp}
In this section, we describe a family of infinite dihedral groups $D_{\infty}$ acting on Cantor sets and give fragmentations of some of them to produce groups whose growth functions are bounded above by $\exp(CR^{\alpha})$ for $\alpha\in(0,1)$. The constructions are based on the prototypical example of the Golden Mean dihedral group described by V.~Nekrashevych in \cite[Section 8]{nekrash18} and \cite[Subsection 5.3.5]{nek22}. 
% See also \cite[Subsection 6.3]{kua24a} for a description of the Golden Mean dihedral group in terms of tile inflations of bounded type.

\subsection{The construction}
Let $\theta\in(0,1)$ be an irrational number. Consider the reflections of the circle $\mathbb{R}/\mathbb{Z}$ with respect to the diameters $\{0,1/2\}$ and $\{\theta/2,(\theta+1)/2\}$:
\[x\cdot a=1-x,  \text{  \quad\ \ \ \quad   }x\cdot b=\theta-x.\]
The composition $ab$ (note we are using right actions) is the rotation $x\mapsto x+\theta$. The orbit of $0$ under the action of $\langle a,b \rangle\cong D_{\infty}$ coincides with that under the action of $\langle ab \rangle\cong \mathbb{Z}$. Let $\mathcal{X}_{\theta}$ be the space obtained by replacing each point $x$ with two points $x-0$ and $x+0$ in the orbit of $0$. Note that $1=0-0$ and $0=0+0$. The set $\mathcal{X}_{\theta}$ is a Cantor set (see \cite[Subsection 1.3.1.2]{nek22}). The action of $\langle a,b\rangle$ naturally lifts to an action of $\mathcal{X}_{\theta}$ by the rules
\[(x-0)\cdot g=x\cdot g+0, \text{  \quad\ \ \ \quad   } (x+0)\cdot g=x\cdot g-0,\] for $g\in\{a,b\}$. We will drop the notations ``$+0$" and ``$-0$" since this will not cause ambiguity. 

Now we encode the points of $\mathcal{X}_{\theta}$ by paths in a Bratteli diagram $\mathsf{B}$. 

Let $\mathcal{A}_0$ be the RK-set\footnote{RK is short for Rokhlin-Kakutani partition of an étale groupoid. We will not present a formal definition here. In this section, an RK-set is a collection of towers and a tower is a collection of disjoint clopen intervals moved by generators of $D_{\infty}$. The union of elements in the towers in the same RK-set form a partition of $\mathcal{X}_{\theta}$. For more detail, see \cite[Subsection 5.2.2]{nek22}.} with towers 
%$\mathcal{F}_{0,0}=\{[0,1-\theta]\}$,  $\mathcal{F}_{1,0}=\{[1-\theta,1]\}$ if $\theta<1/2$, or 
$\mathcal{F}_{0,0}=\{[0,\theta]\}$,  $\mathcal{F}_{1,0}=\{[\theta,1]\}$. Then the level $0$ vertices $V_0:=\{0,1\}$ of $\mathsf{B}$ is identified with $\mathcal{F}_{0,0},\mathcal{F}_{1,0}$, respectively. Let $\mathcal{A}_1$ be the RK-set with towers $\mathcal{F}_{0,1}$, $\mathcal{F}_{1,1}$ constructed as follows. Let $r=\floor{\frac{1}{\theta}}$. Assume, without loss of generality, that $\theta<1/2$, since otherwise we can replace $\theta$ with $1-\theta$. Then $r>1$. The base of $\mathcal{F}_{0,1}$ is $I_0=[0,1-r\theta]$ and the base of $\mathcal{F}_{1,1}$ is $J_0=[1-r\theta,\theta]$. To describe these two towers, let us first recall the construction of the towers $T_{[0,1-r\theta]}$ and $T_{[1-r\theta,\theta]}$ with bases $I_0$ and $J_0$ of the rotation $R_{\theta}:x\mapsto x+\theta$ (see \cite[Subsection 1.3.9]{nek22}). 

The first return map maps $I_0$ to $[(r+1)\theta-1,\theta]$ and $J_0$ to $[0,(r+1)\theta-1]$, and the first return times are $r+1$ and $r$, respectively. Then the towers 
\[T_{[0,1-r\theta]}=\{I_0, I_0\cdot R_{\theta},I_0\cdot R^2_{\theta},I_0\cdot R^3_{\theta},\ldots ,I_0\cdot R^i_{\theta},\ldots ,I_0\cdot R^{r-2}_{\theta},I_0\cdot R^{r-1}_{\theta},I_0\cdot R^{r}_{\theta}\}\]
\begin{align*}
=\{&[0,1-r\theta],[\theta,1-(r-1)\theta],[2\theta,1-(r-2)\theta],[3\theta,1-(r-3)\theta],\ldots ,\\
&[i\theta,1-(r-i)\theta],\ldots ,[(r-2)\theta,1-2\theta],[(r-1)\theta,1-\theta],[r\theta,1]\},
\end{align*}
and 
\[T_{[1-r\theta,\theta]}=\{J_0, J_0\cdot R_{\theta},J_0\cdot R^2_{\theta},I_0\cdot R^3_{\theta},\ldots ,J_0\cdot R^i_{\theta},\ldots ,J_0\cdot R^{r-2}_{\theta},J_0\cdot R^{r-1}_{\theta}\}\]
\begin{align*}
=\{&[1-r\theta,\theta],[1-(r-1)\theta,2\theta,],[1-(r-2)\theta,3\theta],[1-(r-3)\theta,4\theta],\ldots ,\\
&[1-(r-i)\theta,(i+1)\theta],\ldots ,[1-2\theta,(r-1)\theta],[1-\theta,r\theta]\}.
\end{align*}
The sets in $T_{[0,1-r\theta]}$ and $T_{[1-r\theta,\theta]}$ form a partion of $\mathcal{X}_{\theta}$. In the tower $T_{[0,1-r\theta]}$, the transformation $a$ interchanges the intervals $[i\theta,1-(r-i)\theta]$ and $[(r-i)\theta,1-i\theta]$, for $i=0,\ldots ,r$, while $b$ interchanges $[i\theta,1-(r-i)\theta]$ and $[(r+1-i)\theta,1-(i-1)\theta]$, for $i=1,\ldots ,r$. In the tower $T_{[1-r\theta,\theta]}$, the transformation $a$ interchanges the intervals $[1-(r-i)\theta,(i+1)\theta]$ and $[1-(i+1)\theta,(r-i)\theta]$, for $i=0,\ldots ,r-1$, while $b$ interchanges $[1-(r-i)\theta,(i+1)\theta]$ and $[1-i\theta,(r-i+1)\theta]$, for $i=1,\ldots r-1$. It follows that the towers $\mathcal{F}_{0,1}$ and
$\mathcal{F}_{1,1}$ consist of the same elements as those of $T_{[0,1-r\theta]}$ and $T_{[1-r\theta,\theta]}$, respectively. The arrows between the intervals are shown in Figure \ref{rksets}. Note that the intervals are lined up in the order of iterations of $R_{\theta}$ before the first return. 

The level $1$ vertices $V_1$ of $\mathsf{B}$ consists of $2$ elements $0,1$ corresponding to $\mathcal{F}_{0,1}$ and $\mathcal{F}_{1,1}$, respectively. The base $[0,\theta]\in\mathcal{F}_{0,0}$ is partitioned by $I_0,J_0$ while the base $[\theta,1]\in \mathcal{F}_{1,0}$ is partitioned by the rest of elements in $\mathcal{F}_{0,1}$ and $\mathcal{F}_{1,1}$. Hence the finite (unordered) diagram of $\mathsf{B}$ associated with the sequence $\mathcal{A}_0\prec \mathcal{A}_1$ consists of vertices $V_0,V_1$ where $0\in V_0$ is connected to $0,1\in V_1$ each by a single edge, and $1\in V_0$ is connected to $0\in V_1$ by $r$ edges and to $1\in V_1$ by $r-1$ edges. See Figure \ref{bratlv1} for an example when $r=5$. 

\begin{figure}

\centering
\begin{minipage}{1\textwidth}
\centering
\subcaption{Tower $\mathcal{F}_{0,1}$.}
\begin{adjustbox}{width=.55\textwidth}
\begin{tikzpicture}
\begin{scope}state without output/.append style={draw=none}%[every node/.style={circle,thick,draw}]
    \node (A) at (0,0) {$[0,1-r\theta]$};
    \node (B) at (0,-1) {$[\theta,1-(r-1)\theta]$};
    \node (C) at (0,-2) {$[2\theta,(1-r-2)\theta]$};
    \node (E) at (0,-3) {$[3\theta,1-(r-3)\theta]$};
    \node (F) at (0,-4) {$[4\theta,1-(r-4)\theta]$};
    \node (G) at (0,-6) {\ldots };
    \node (Z) at (0,-8) {\ldots }; 
    \node (H) at (0,-9) {$[(r-3)\theta,1-3\theta]$};
    \node (I) at (0,-10) {$[(r-2)\theta,1-2\theta]$};
    \node (J) at (0,-11) {$[(r-1)\theta,1-\theta]$};
    \node (K) at (0,-12) {$[r\theta,1]$}; 
\end{scope}

\begin{scope} %[>={Stealth[black]},
              every node/.style={fill=white,circle},
              every edge/.style={draw=black,very thick}]
              
    %\path [->] (A) edge [bend left=40] node {} (B);
    %\path [->] (B) edge [bend left=40] node {} (A);
    %\draw [->] (B) to  [out=320,in=40,looseness=9] (B);
    %\path [->] (D) edge node {$3$} (C);
    %\path [->] (A) edge node {$3$} (E);
    %\path [->] (D) edge node {$3$} (E);
    %\path [->] (D) edge node {$3$} (F);
    %\path [->] (C) edge node {$5$} (F);
    %\path [->] (E) edge node {$8$} (F); 
    %\path [->] (B) edge[bend right=60] node {$1$} (E); 

    %\draw[]  (A) node[draw=none][midway,above] {$b_1$} (B);
    \path [<->] (A) edge [bend right=99] node  [midway,left] {$a$} (K);
    \path [<->] (B) edge [bend right=80] node [midway,right] {$a$} (J);
    \path [<->] (C) edge [bend right=75] node  [midway,left] {$a$} (I);
    \path [<->] (E) edge [bend right=65] node  [midway,left] {$a$} (H);
    \path [<->] (F) edge [bend right=60] node  [midway,left] {$a$}  (Z);
   
    \path [<->] (B) edge [bend left=99] node  [midway,right] {$b$} (K);
    \path [<->] (C) edge [bend left=80] node  [midway,left] {$b$} (J);
    \path [<->] (E) edge [bend left=75] node  [midway,right] {$b$} (I);
    \path [<->] (F) edge [bend left=65] node  [midway,right] {$b$} (H);
    
\end{scope}

\end{tikzpicture}
\end{adjustbox}

\label{tower1}
%\end{figure}
\end{minipage}

%\begin{figure}

\begin{minipage}{1\textwidth}
\centering
\subcaption{Tower $\mathcal{F}_{1,1}$.}
\begin{adjustbox}{width=.55\textwidth}
\begin{tikzpicture}
\centering
\begin{scope}state without output/.append style={draw=none}%[every node/.style={circle,thick,draw}]
    \node (A) at (0,0) {$[1-r\theta,\theta]$};
    \node (B) at (0,-1) {$[1-(r-1)\theta,2\theta]$};
    \node (C) at (0,-2) {$[1-(r-2)\theta,3\theta]$};
    \node (E) at (0,-3) {$[1-(r-3)\theta,4\theta]$};
    \node (F) at (0,-4) {$[1-(r-4)\theta,5\theta]$};
    \node (G) at (0,-6) {\ldots };
    \node (Z) at (0,-8) {\ldots }; 
    \node (H) at (0,-9) {$[1-4\theta,(r-3)\theta]$};
    \node (I) at (0,-10) {$[1-3\theta,(r-2)\theta]$};
    \node (J) at (0,-11) {$[1-2\theta,(r-1)\theta]$};
    \node (K) at (0,-12) {$[1-\theta,r\theta]$}; 
\end{scope}

\begin{scope} %[>={Stealth[black]},
              every node/.style={fill=white,circle},
              every edge/.style={draw=black,very thick}]
              
    %\path [->] (A) edge [bend left=40] node {} (B);
    %\path [->] (B) edge [bend left=40] node {} (A);
    %\draw [->] (B) to  [out=320,in=40,looseness=9] (B);
    %\path [->] (D) edge node {$3$} (C);
    %\path [->] (A) edge node {$3$} (E);
    %\path [->] (D) edge node {$3$} (E);
    %\path [->] (D) edge node {$3$} (F);
    %\path [->] (C) edge node {$5$} (F);
    %\path [->] (E) edge node {$8$} (F); 
    %\path [->] (B) edge[bend right=60] node {$1$} (E); 

    %\draw[]  (A) node[draw=none][midway,above] {$b_1$} (B);
    \path [<->] (A) edge [bend right=99] node  [midway,left] {$a$} (K);
    \path [<->] (B) edge [bend right=80] node [midway,right] {$a$} (J);
    \path [<->] (C) edge [bend right=75] node  [midway,left] {$a$} (I);
    \path [<->] (E) edge [bend right=65] node  [midway,left] {$a$} (H);
    \path [<->] (F) edge [bend right=60] node  [midway,left] {$a$}  (Z);
   
    \path [<->] (B) edge [bend left=99] node  [midway,right] {$b$} (K);
    \path [<->] (C) edge [bend left=80] node  [midway,left] {$b$} (J);
    \path [<->] (E) edge [bend left=75] node  [midway,right] {$b$} (I);
    \path [<->] (F) edge [bend left=65] node  [midway,right] {$b$} (H);
    
\end{scope}

\end{tikzpicture}
\end{adjustbox}

\label{tower2}
\end{minipage}
\caption{RK-sets of dihedral groups}
\label{rksets}
\end{figure}

\begin{figure}[!htb]
\centering
\begin{tikzpicture}
\begin{scope}state without output/.append style={draw=none}%[every node/.style={circle,thick,draw}]
    \node (A) at (0,0) {};
    \node (B) at (5,0) {};
    \node (C) at (0,-6) {};
    \node (D) at (5,-6) {};

\end{scope}

%\begin{scope} %[>={Stealth[black]},
              every node/.style={fill=white,circle},
              every edge/.style={draw=black,very thick}]
              
    %\path [->] (A) edge [bend left=40] node {} (B);
    %\path [->] (B) edge [bend left=40] node {} (A);
    %\draw [->] (B) to  [out=320,in=40,looseness=9] (B);
    %\path [->] (D) edge node {$3$} (C);
    %\path [->] (A) edge node {$3$} (E);
    %\path [->] (D) edge node {$3$} (E);
    %\path [->] (D) edge node {$3$} (F);
    %\path [->] (C) edge node {$5$} (F);
    %\path [->] (E) edge node {$8$} (F); 
    %\path [->] (B) edge[bend right=60] node {$1$} (E); 

    %\draw[]  (A) node[draw=none][midway,above] {$b_1$} (B);

%\end{scope}

\begin{scope}
\path [-] (A) edge node [midway,left] {} (C);
\path [-] (B) edge [bend right=20] node [near end,left] {} (C);
\path [-] (B) edge [bend right=10] node [near end,left] {} (C);
\path [-] (B) edge node [near end,left] {} (C);
\path [-] (B) edge [bend left=10] node [near end,left] {} (C);
\path [-] (B) edge [bend left=20] node [near end,left] {} (C);
\path [-] (B) edge [bend left=15] node [near end,left] {} (D);
\path [-] (B) edge [bend left=5] node [near end,left] {} (D);
\path [-] (B) edge [bend right=5] node [near end,left] {} (D);
\path [-] (B) edge [bend right=15] node [near end,left] {} (D);
\path [-] (D) edge node [near end,right]{} (A);
% \draw (0,-5) node {} -- (0,-6) [dotted] node {};
 % \draw (4,-5) node {} -- (4,-6) [dotted] node {};

%\path [-] (B) edge [bend right=35] node  [midway,left]{$1$} (D) ;
%\path [-] (B) edge [bend right=20] node [midway,left]{$2$} (D);
%\path [-] (B) edge [bend right=5] node [midway,left]{$3$} (D);
%\path [-] (B) edge [bend left=5] node [midway,right]{$4$} (D);
%\path [-] (B) edge [bend left=20] node [midway,right]{$5$} (D);
%\path [-] (B) edge [bend left=35] node [midway,right]{$6$} (D);

\filldraw[black] (0,0) circle  (2.5pt) node[anchor=west] {0};
\filldraw[black] (5,0) circle (2.5pt) node[anchor=west] {1};
\filldraw[black] (0,-6) circle (2.5pt) node[anchor=west] {0};
\filldraw[black] (5,-6) circle (2.5pt) node[anchor=west] {1};
\end{scope}

\end{tikzpicture}
\caption{Finite diagram for $\mathcal{A}_0\prec\mathcal{A}_1$ when $r=5$.}
\label{bratlv1}
\end{figure}
Now consider the similarity map $S:[0,\theta]\rightarrow [0,1]$ given by $x\mapsto x/\theta$. Then $S([0,1-r\theta])=[0,1/\theta-r]$ and $S([1-r\theta,\theta]=[1/\theta-r,1])$, and they form a partition of $[0,1]$. Denote $\theta'=1/\theta-r$. The actions of $a,b$ on $[0,\theta]$ are conjugated by $S$ to be the actions $x\mapsto \theta'-x$ and $x\mapsto 1-x$, respectively, for $x\in\mathcal{X}_{\theta'}$. Indeed, this follows from the following proposition. 

\begin{prop} We have the following: 
\begin{enumerate}
    \item The actions of $a$ on $[0,1-r\theta]$ and $[1-r\theta,\theta]$ are conjugated by $S$ to be $1/\theta-r-x$ for $x\in[0,1/\theta-r]$ and $[1/\theta-r,1]$, respectively.
    \item The actions of $b$ on $[0,1-r\theta]$ and $[1-r\theta,\theta]$ are conjugated by $S$ to be $1-x$ for $x\in[0,1/\theta-r]$ and $[1/\theta-r,1]$, respectively.
\end{enumerate}
    
\end{prop}
\begin{proof}
    The above statements follow from the following $4$ commutative diagrams.
\[\begin{tikzcd}
\left[0,1-r\theta\right] \arrow{r}{a} \arrow[swap]{d}{S} & \left[r\theta,1\right] \arrow{d}{S} \\%
\left[0,1/\theta-r\right] \arrow{r}{a}& \left[r,1/\theta\right]=\left[0,1/\theta-r\right]
\end{tikzcd}
\]
\vspace{.25cm}
\[ \begin{tikzcd}
\left[1-r\theta,\theta\right] \arrow{r}{a} \arrow[swap]{d}{S} & \left[1-\theta,r\theta\right] \arrow{d}{S} \\%
\left[1/\theta-r,1\right] \arrow{r}{a}& \left[1/\theta-1,r\right]=\left[1/\theta-r,1\right]
\end{tikzcd}
\]
\vspace{.5cm}
\[ \begin{tikzcd}
\left[0,1-r\theta\right] \arrow{r}{b} \arrow[swap]{d}{S} & \left[(1+r)\theta-1,\theta\right] \arrow{d}{S} \\%
\left[0,1/\theta-r\right] \arrow{r}{b}& \left[1-(1/\theta-r),1\right]
\end{tikzcd}
\]
\vspace{.25cm}
\[ \begin{tikzcd}
\left[1-r\theta,\theta\right] \arrow{r}{b} \arrow[swap]{d}{S} & \left[0,(1+r)\theta-1\right] \arrow{d}{S} \\%
\left[1/\theta-r,1\right] \arrow{r}{b}& \left[0,1-(1/\theta-r)\right]
\end{tikzcd}
\]

Since $a,b$ have order $2$, $a$ flips the intervals on the bottom row of the first $2$ diagrams, and $b$ interchanges the intervals on the bottom row of the last $2$ diagrams. Hence $a$ acts as $1/\theta-r-x=\theta'-x$ and $b$ acts as $1-x$ for $x\in S([0,\theta])=[0,1]$. 
\end{proof}

Let $r'=\floor{\frac{1}{\theta'}}$. Let $I_0'=[0,1-r'\theta']$ and $J_0'=[1-r'\theta,\theta']$. Let $\mathcal{A}_2$ be the RK-set associated with the new $a,b$ with $2$ towers $\mathcal{F}_{0,2}$ and $\mathcal{F}_{1,2}$ whose bases are $I_0'$ and $J_0'$, respectively. Then the finite diagram of $\mathsf{B}$ associated with $\mathcal{A}_1\prec\mathcal{A}_2$ is determined by the same procedure as that of $\mathcal{A}_0\prec\mathcal{A}_1$, with the roles of $a$ and $b$ interchanged in Figure \ref{rksets}. 

Let $\theta_0=\theta$, $\theta_1=\theta'$, $c_0=r$, and $c_1=r'$. Then we have \[\theta_0=\dfrac{1}{c_0+\theta_1}.\] It follows that the diagram associated with $\mathcal{A}_0\prec\mathcal{A}_1$ is determined by the procedure above for the maps $x\mapsto 1-x$ and $x\mapsto \theta_0-x$. In general, the finite diagram associated with $\mathcal{A}_{n-1}\prec\mathcal{A}_{n}$ is determined by the same procedure for the maps $x\mapsto 1-x$ and $x\mapsto \theta_{n-1}-x$ where $\theta_{n-1}$ is defined by the rule \[\theta_{n-1}=\dfrac{1}{c_{n-1}+\theta_n}, \text{ for }n\in\mathbb{N},\] where $c_{n-1}=\floor{\frac{1}{\theta_{n-1}}}$. For the RK-set $\mathcal{A}_{n}$ ($n\geq 1$), its bases of towers are denoted $I_0^{(n)}=[0,1-c_{n-1}\theta_{n-1}]$ and $J_0^{(n)}=[1-c_{n-1}\theta_{n-1},\theta_{n-1}]$. Therefore, we have the following. 

\begin{thm}\label{makebratt}
    Let $\theta\in(0,1)$ be an irrational number, and let $(c_0,c_1,c_2,\ldots )$ be the sequence of positive integers such that 
    \[\theta=\dfrac{1}{c_0+\dfrac{1}{c_1+\dfrac{1}{c_2+\dfrac{1}{\ldots }}}}.\]

 The (unordered) Bratteli diagram $\mathsf{B}$ associated with the sequence of RK-sets $\mathcal{A}_0\prec\mathcal{A}_1\prec\mathcal{A}_2\prec\ldots $ is described as follows. The $n$-th level vertices $V_n=\{0,1\}$ for $n\geq 0$. If $\theta_{n-1}<1/2$ (equivalently $c_{n-1}\geq 2$), then $0\in V_{n-1}$ is connected to $0,1\in V_n$ each by a single edge, and $1\in V_{n-1}$ is connected to $0\in V_n$ by $c_{n-1}$ edges and to $1\in V_n$ by $c_{n-1}-1$ edges. If $\theta_{n-1}>1/2$ (equivalently $c_{n-1}=1$), then $0\in V_{n-1}$ is connected to $0,1\in V_{n}$ each by a single edge, and $1\in V_{n-1}$ is connected to $0\in V_{n}$ by a single edge. The set of edges on each level of $\mathsf{B}$ is denoted $E_n$, for $n=1,2,\ldots $, satisfying $\mathbf{s}(E_n)=V_{n-1}$ and $\mathbf{r}(E_n)=V_{n}$. 

 The natural action $\Omega(\mathsf{B})\curvearrowleft D_{\infty}$ is topologically conjugate to the action $\mathcal{X}_{\theta}\curvearrowleft D_{\infty}$. 
\end{thm}
\begin{proof}
  The description of $\mathsf{B}$ follows from the construction above and the definition of Bratteli diagrams associated with a sequence of RK-sets. Let us prove that $\Omega(\mathsf{B})\curvearrowleft D_{\infty}$ is topologically conjugate to the action $\mathcal{X}_{\theta}\curvearrowleft D_{\infty}$. The elements of each $E_n$ will be denoted $I_k^{(n)}$ and $J_l^{(n)}$, using the same notations as elements in $\mathcal{F}_{0,n}$ and $\mathcal{F}_{1,n}$, respectively, for $c_{n-1}>2$, $0\leq k\leq c_{n-1}$ and $0\leq l\leq c_{n-1}-1$. If $c_{n-1}=1$, then $E_n=\{I_0^{(n)},J_0^{(n)},I_1^{(n)}\}$. The lower indices $k,l$ appear in order of successive applications of $abab\ldots $.  %\textcolor{red}{In both cases, $I_1^{(n)}=I_0^{(n)}\cdot a$ if $n$ is odd, and $I_0^{(n)}\cdot b$ if $n$ is even.(DELETE need symbols for each $\theta_{n-1}-x$ and $1-x$)}

\begin{lemma}\label{conjugator}
    Let $\mathsf{s}:\Omega(\mathsf{B})\rightarrow \mathsf{s}(\Omega(\mathsf{B}))$ be the shift map, i.e., $\mathsf{s}(\ldots e_2e_1)=\ldots e_3e_2$. Let $\{\lambda_i:\mathsf{s}^{i-1}(\Omega(\mathsf{B}))\rightarrow\mathcal{X}_{\theta_{i-1}}\}_{i\in\mathbb{N}}$ be a family of maps given by the condition $\lambda_i(\ldots e_{i+1}e_i)\in e_i$ together with the recurrent rule
    \begin{align}
   \lambda_{i+1}(\ldots e_{i+2}e_{i+1})= 
   \begin{cases}
       \frac{1}{\theta_{i-1}}\lambda_{i}(\ldots e_{i+1}e_i), &\text{ if } e_i\in\{I_{0}^{(i)}, J_{0}^{(i)}\},\\
      \frac{1}{\theta_{i-1}}(1-\lambda_i(\ldots e_{i+1}e_i)\cdot g_{e_i}), &\text{ if }e_i\in \mathcal{F}_{0,i}\backslash\{I_{0}^{(i)}\},\\
      \frac{1}{\theta_{i-1}}(1-\lambda_i(\ldots e_{i+1}e_i)\cdot h_{e_i}), &\text{ if }e_i\in \mathcal{F}_{1,i}\backslash\{J_{0}^{(i)}\},\\
   \end{cases}
\end{align}
where $g_{e_i}\in D_{\infty}$ is such that $e_i\cdot g_{e_i}=I_{1}^{(i)}$ with $|g_{e_i}|<c_{i-1}$, and $h_{e_i}\in D_{\infty}$ is such that $e_i\cdot h_{e_i}=J_{1}^{(i)}$ with $|h_{e_i}|<c_{i-1}-1$, if $c_{i-1}\geq 2$. If $c_{i-1}=1$, then $g_{e_i}=Id$ and the third row of the recursion must be deleted since in this case $\mathcal{F}_{1,i}=\{J_{0}^{(i)}\}$. Then each $\lambda_i$ is a homeomorphism. 
\end{lemma}
\begin{proof}
Let us first check that each $\lambda_i$ is well-defined, i.e., for each $\ldots e_2e_1\in\Omega(\mathsf{B})$, the value $\lambda_1(\ldots e_2e_1)$ is unique. Note that each $e_i$ also represents a component (an interval) of a tower in $\mathcal{X}_{\theta_{i-1}}$, and thus $e_i$ is also treated as a set. By the definition, $\lambda_i(\ldots e_{i+1}e_{i})\in e_i$, $\lambda_{i+1}(\ldots e_{i+2}e_{i+1})\in e_{i+1}$, and  
 \begin{align}\label{recrr}
   \lambda_{i}(\ldots e_{i+1}e_{i})= 
   \begin{cases}
       \theta_{i-1}\lambda_{i+1}(\ldots e_{i+2}e_{i+1}), &\text{ if } e_i\in\{I_{0}^{(i)}, J_{0}^{(i)}\},\\
     (1-\theta_{i-1}\lambda_{i+1}(\ldots e_{i+2}e_{i+1}))\cdot g_{e_i}^{-1}, &\text{ if }e_i\in \mathcal{F}_{0,i}\backslash\{I_{0}^{(i)}\},\\
       (1-\theta_{i-1}\lambda_{i+1}(\ldots e_{i+2}e_{i+1}))\cdot h_{e_i}^{-1}, &\text{ if }e_i\in \mathcal{F}_{1,i}\backslash\{J_{0}^{(i)}\}.\\
   \end{cases}
\end{align}
It follows that $e_{i+1}$ is mapped by the rule $(\ast)$ to be a subset of $e_i$ where $(\ast)$ is obtained by replacing $\lambda_{i+1}(\ldots e_{i+2}e_{i+1})$ with $e_{i+1}$ on the right hand side of (\ref{recrr}). Let $C_1=e_1$, $C_2$ be obtained from applying $(\ast)$ on $e_2$ once, and in general $C_i$ be obtained from applying $(\ast)$ on $e_i$ for $i-1$ times. Then we constructed a sequence of clopen sets $C_1\supset C_2\supset\ldots $ whose diameters converge to $0$. Since each $C_i$ is compact, $\bigcap_{i=1}^{\infty}C_i$ consists of a single point. It follows that the map $\lambda_1$ (as well as $\lambda_i$ for all $i\in \mathbb{N}$) is well defined. Since the elements in the RK-set $\mathcal{A}_{i-1}$ forms a partition of $\mathcal{X}_{\theta_{i-1}}$, for $i\in\mathbb{N}$, each point $x\in \mathcal{X}_{\theta_{i-1}}$ is uniquely determined by its address $(C_i,C_{i+1},\ldots )$. It follows that the maps $\lambda_i$ are both injective and surjective. The map $\lambda_1$ being continuous follows from the fact that $x_1,x_2\in\bigcap_{i=1}^n C_i$ if and only if their preimages under $\lambda_1$, denoted $\gamma_1,\gamma_2$, share common prefixes up to at least level $n$ of $\mathsf{B}$. The proof is complete. 
\end{proof}
 
\begin{lemma}\label{natuac}
 Suppose all $c_i\geq 2$. For each $\mathcal{X}_{\theta_{i-1}}$, denote by $a^{(i)}=1-x$ and $b^{(i)}=\theta_{i-1}-x$, for $x\in\mathcal{X}_{\theta_{i-1}}$ and $i\geq 1$. The natural action $\Omega(\mathsf{B})\curvearrowleft D_{\infty}$ is given by the following rules. 
  \begin{enumerate}
      \item\label{ruleinter} Let $e\in E_i$. If $e\in\{I_0^{(i)},J_0^{(i)}\}$, then ${}_{e}{|a^{(i)}}=b^{(i+1)}$ and ${}_{e}{|b^{(i)}}=a^{(i+1)}$. If $c_{i-1}$ is odd, then ${}_{I_{c_{i-1}}^{(i)}}{|b^{(i)}}=b^{(i+1)}$ and ${}_{J_{c_{i-1}-1}^{(i)}}{|a^{(i)}}=b^{(i+1)}$. If $c_{i-1}$ is even, then ${}_{I_{c_{i-1}}^{(i)}}{|a^{(i)}}=b^{(i+1)}$ and ${}_{J_{c_{i-1}-1}^{(i)}}{|b^{(i)}}=b^{(i+1)}$. 
    \item\label{deter1} Suppose $c_0$ is odd. For $w\in\mathsf{s}(\Omega(\mathsf{B}))$, $a_1$ interchanges the elements $wI_{2k}^{(1)}$ with $wI_{2k+1}^{(1)}$, for $k=0,\ldots ,(c_0-1)/2$, and interchanges $wJ_{2k}^{(1)}$ with $wJ_{2k+1}^{(1)}$, for $k=0,\ldots ,(c_0-3)/2$; $b_1$ interchanges the elements $wI_{2k+1}^{(1)}$ with $wI_{2k+2}^{(1)}$, for $k=0,\ldots ,(c_0-3)/2$, and interchanges $wJ_{2k+1}^{(1)}$ with $wJ_{2k+2}^{(1)}$, for $k=0,\ldots ,(c_0-3)/2$. Their sections on the first edges are identity.  %Moreover, $a_1$ fixes $J_{c_0-1}^{(1)}$ and $b_1$ fixes $I_{c_0}^{(1)}$, with non-identity sections to be described below. 
    
    \item\label{deter2} Suppose $c_0$ is even. For $w\in\mathsf{s}(\Omega(\mathsf{B}))$, $a_1$ interchanges the elements $wI_{2k}^{(1)}$ with $wI_{2k+1}^{(1)}$, for $k=0,\ldots ,c_0/2$, and interchanges $wJ_{2k}^{(1)}$ with $wJ_{2k+1}^{(1)}$, for $k=0,\ldots ,c_0/2$; $b_1$ interchanges the elements $wI_{2k+1}^{(1)}$ with $wI_{2k+2}^{(1)}$, for $k=0,\ldots ,c_0/2$, and interchanges $wJ_{2k+1}^{(1)}$ with $wJ_{2k+2}^{(1)}$, for $k=0,\ldots ,(c_0-2)/2$. Their sections on the first edges are identity.  %Moreover, $a_1$ fixes $I_{c_0}^{(1)}$ and $b_1$ fixes $J_{c_0-1}^{(1)}$, with non-identity sections to be described below. 

\item\label{nondeter1} For $w\in\mathsf{s}^{i+2}(\Omega(\mathsf{B}))$,
  \begin{align*}
 & wI_0^{(i+1)}I_0^{(i)}\cdot b^{(i)} = wI_0^{(i+1)}\cdot a^{(i+1)}J_0^{(i)}=wI_1^{(i+1)}J_0^{(i)},\\
  & wI_1^{(i+1)}J_0^{(i)}\cdot b^{(i)}= wI_1^{(i+1)}\cdot a^{(i+1)}I_0^{i}= wI_0^{(i+1)}I_0^{(i)},\\
  & wJ_0^{(i+1)}I_0^{(i)}\cdot b^{(i)}=wJ_0^{(i+1)}\cdot a^{(i+1)}J_{0}^{i}=wJ_1^{(i+1)}J_0^{(i)},\\
  & wJ_1^{(i+1)}J_0^{(i)}\cdot b^{(i)}=wJ_1^{(i+1)}\cdot a^{(i+1)} I_0^{(i)}=wJ_0^{(i+1)}I_0^{(i)}.\\
\end{align*}
\item\label{nonodeter2} For $w\in\mathsf{s}^{i+2}(\Omega(\mathsf{B}))$, if $c_{i}$ is odd, then
\[wI_{c_i}^{(i+1)}J_0^{(i)}\cdot b^{(i)}=wI_{c_i}^{(i+1)}\cdot a^{(i+1)}J_0^{(i)}=w\cdot b^{(i+2)}I_{c_i}^{(i+1)}J_0^{(i)}.\]
If $c_{i}$ is even, then
\[wJ_{c_{i}-1}^{(i+1)}J_0^{(i)}\cdot b^{(i)}=wJ_{c_i-1}^{(i+1)}\cdot a^{(i+1)}J_0^{(i)}=w\cdot b^{(i+2)}J_{c_i-1}^{(i+1)}J_0^{(i)}.\]
\end{enumerate}

\end{lemma}   
   % For $w\in\mathsf{s}(\Omega(\mathsf{B}))$ and $c_0>2$,
%
  %  w1\cdot a &= (w\cdot b)1,\text{      } w12\cdot a= (w\cdot b)12,\\
  %  w11\cdot b =  w12,&\text{   } w12\cdot b=w11,\text{   } w121\cdot b=(w\cdot b)121. 
Rules \ref{nondeter1} and \ref{nonodeter2} show that each $b^{(i)}$ is nondeterministic at each $J_0^{(i)}$. 
\begin{proof}
    Notice that $a_1=a$ and $b_1=b$. The first sentence in Rule \ref{ruleinter} follows directly from the fact that the similarity map $S$ interchanges the roles of $a^{(i)},b^{(i)}$ on $I_0^{(i)}\cup J_0^{(i)}$, i.e., $S\circ a^{(i)}=\theta_i-x=b^{(i+1)}$ and $S\circ b^{(i)}=1-x=a^{(i+1)}$. Moreover, if $c_0$ is odd the interval $I_{c_0}^{(i)}$ is flipped by $b^{(i)}$, and the interval $J_{c_0-1}^{(i)}$ is flipped by $a^{(i)}$, while they are mapped by $\lambda_i$ onto $[0,\theta_i]$ and $[\theta_i,1]$, respectively. Hence both $b^{(i)}$ and $a^{(i)}$ are conjugated to be $b^{(i+1)}$ on $[0,\theta_i]$ and $[\theta_i,1]$, respectively. If $c_0$ is even, the argument is similar by noting that the interval $I_{c_0}^{(i)}$ is flipped by $a^{(i)}$, and the interval $J_{c_0-1}^{(i)}$ is flipped by $b^{(i)}$. Rules \ref{deter1} and \ref{deter2} follow directly from the description of the RK-sets shown in Figure \ref{rksets}. 

    For Rules \ref{nondeter1} and \ref{nonodeter2}, notice that we assume that all $c_i\geq2$, which is equivalent to all $\theta_i<1/2$. Hence $I_0^{(i)}\cdot b^{(i)}=[0,1-c_{i-1}\theta_{i-1}]\cdot b^{(i)}=[(c_{i-1}+1)\theta_{i-1}-1,\theta_{i-1}]\subset J_0^{(i)}$. It follows that $b^{(i)}$ interchanges the intervals $I_0^{(i)}$ with $[(c_{i-1}+1)\theta_{i-1}-1,\theta_{i-1}]$ and flips the interval $[1-c_{i-1}\theta_{i-1},(c_{i-1}+1)\theta_{i-1}-1]\subset J_0^{(i)}$. It follow that on the corresponding part of the Bratteli diagram, $b^{(i)}$ has $2$ possibilities on the edge $J_0^{(i)}$: (1) mapping it to $I_0^{(i)}$; (2) fixing it. The sections of $b^{(i)}$ on the corresponding edges follow from Rule \ref{ruleinter}.  
\end{proof}
By the above lemmas, the natural action $\Omega(\mathsf{B})\curvearrowleft D_{\infty}$ is conjugated to $\mathcal{X}_{\theta_0}\curvearrowleft D_{\infty}$ by the map $\lambda_1$. The case when some $c_i=1$ was shown in \cite[Subsection 5.3.4]{nek22}. We omit it here. The proof is complete. 
\end{proof}

\subsection{Actions $\Omega(\mathsf{B})\curvearrowleft D_{\infty}$ in terms of tile inflations}
Let $(c_0,c_1,\ldots )$ be a sequence of positive integers as in Theorem \ref{makebratt}. As before, we assume that each $c_i\geq 2$. In this subsection, we convert the action described in Lemma \ref{natuac} to tile inflations. Let $\mathcal{A}$ be the automaton whose initial states generate $D_{\infty}$ as an inverse semigroup. In other words, the initial states of $\mathcal{A}$ are $a,b$ restricted at each edge in $E_1$. We will not split the generators $a,b$ (as well as each $a^{(i)},b^{(i)}$), as it is clear that when restricting them at a single edge, the resulting partial transformation is a state in $\mathcal{A}$. 

Since each $V_n=\{0,1\}$, there are $2$ tiles on each level of $\mathsf{B}$, denoted $\mathcal{T}_{0,n}$ and $\mathcal{T}_{1,n}$. The tile $\mathcal{T}_{0,1}$ is linear whose vertices and edges are shown in Figure \ref{tower1} in the order of successive applications of $abab\ldots $ starting from $I_0^{(1)}=[0,1-r\theta]$. The tile $\mathcal{T}_{1,1}$ is also linear whose vertices and edges are shown in Figure \ref{tower2} in the order of successive applications of $abab\ldots $ starting from $J_0^{(1)}=[1-r\theta,\theta]$. The generator $a$ has $1$ fixed point and $b$ has $2$ fixed points. To describe the fixed points (boundary points on infinite tiles), boundary points on finite tiles and connecting points, we need to distinguish $4$ cases according to the parities of adjacent $c_{i},c_{i+1}$. For convenience, we only describe in detail Case 1 below for $c_0,c_1$. The other $3$ cases (together with Case 1) are summarized in Table \ref{allcase1}. The general descriptions for $c_{i},c_{i+1}$ will follow from the same patterns. Especially, the continuations of boundary points and locations of connecting points are obtained by successively reading off the tables. 
%\begin{casess}[leftmargin=2.7\parindent] 

    \textbf{Case 1}:\label{oddeven} $c_0$ is odd and $c_1$ is even. On $E_1$, $a$ fixes $J_{c_0-1}^{(1)}$ and $b$ fixes $I_{c_0}^{(1)}$, each has nontrivial sections. Moreover, $a$ has nontrivial section at $I_0^{(1)}$ and $b$ has nontrivial section at $J_0^{(1)}$ (recall that $b$ is nondeterministic at $J_0^{(1)}$). Therefore, the boundary points of $\mathcal{T}_{0,1}$ are $I_0^{(1)}$, $I_{c_0}^{(1)}$, and the boundary points of $\mathcal{T}_{0,1}$ are $J_0^{(1)}$, $J_{c_0-1}^{(1)}$. The tile $\mathcal{T}_{0,2}$ is obtained as follows. Take $1$ copy $I_0^{(2)}\mathcal{T}_{0,1}$ of $\mathcal{T}_{0,1}$ and $c_1$ copies of $\mathcal{T}_{1,1}$ appending respectively $I_k^{(1)}$, for $k=1,\ldots ,c_1$, denoted $I_k^{(2)}\mathcal{T}_{1,1}$. Connect $I_0^{(2)}I_0^{(1)}$ on $I_0^{(2)}\mathcal{T}_{0,1}$ to $I_1^{(2)}J_0^{(1)}$ on $I_1^{(2)}\mathcal{T}_{1,1}$ by a two-sided arrow labeled by $b$. Connect $I_{2k+1}^{(2)}J_{c_0-1}^{(1)}$ on $I_{2k+1}^{(2)}\mathcal{T}_{1,1}$ to $I_{2k+2}^{(2)}J_{c_0-1}^{(1)}$ on $I_{2k+2}^{(2)}\mathcal{T}_{1,1}$ by a two-sided arrow labeled by $a$, for $k=0,\ldots ,c_1/2$. Connect $I_{2k}^{(2)}J_{0}^{(1)}$ on $I_{2k}^{(2)}\mathcal{T}_{1,1}$ to $I_{2k+1}^{(2)}J_{0}^{(1)}$ on $I_{2k+1}^{(2)}\mathcal{T}_{1,1}$ by a two-sided arrow labeled by $b$, for $k=1,\ldots ,c_1/2$. The boundary points of $\mathcal{T}_{0,2}$ are $I_0^{(2)}I_{c_0}^{(1)}$ and $I_{c_1}^{(2)}J_{0}^{(1)}$. Similarly, the tile $\mathcal{T}_{1,2}$ is obtained as follows. Take $1$ copy $J_0^{(2)}\mathcal{T}_{0,1}$ of $\mathcal{T}_{0,1}$ and $c_1-1$ copies of $\mathcal{T}_{1,1}$, appending respectively $J_{k}^{(2)}$, for $k=1,\ldots ,c_1-1$, denoted $J_{k}^{(2)}\mathcal{T}_{1,1}$. Connect $J_0^{(2)}I_0^{(1)}$ on $\mathcal{T}_{0,1}$ to $J_1^{(2)}J_0^{(1)}$ on $J_{1}^{(2)}\mathcal{T}_{1,1}$ by a two-sided arrow labeled by $b$. Connect $J_{2k+1}^{(2)}I_{c_0-1}^{(1)}$ on $J_{2k+1}^{(2)}\mathcal{T}_{1,1}$ to $J_{2k+2}^{(2)}I_{c_0-1}^{(1)}$ on $J_{2k+2}^{(2)}\mathcal{T}_{1,1}$ for $k=0,\ldots ,(c_1-2)/2$ by a two-sided arrow labeled by $a$. Connect $J_{2k}^{(2)}J_0^{(1)}$ on $J_{2k}^{(2)}\mathcal{T}_{1,1}$ to $J_{2k+1}^{(2)}J_0^{(1)}$ on $J_{2k+1}^{(2)}\mathcal{T}_{1,1}$ by a two-sided arrow labeled by $b$, for $k=1,\ldots ,(c_1-2)/2$. The boundary points of $\mathcal{T}_{1,2}$ are $J_0^{(2)}I_{c_1}^{(1)}$ and $J_{c_1-1}^{(2)}J_{c_0-1}^{(1)}$. See Table \ref{case1}. 
    
    \textbf{Case 2}: $c_0$ is even and $c_1$ is odd. See Table \ref{case2}. 
    
    \textbf{Case 3}: both $c_0$ and $c_1$ are odd. See Table \ref{case3}.
    
    \textbf{Case 4}: both $c_0$ and $c_1$ are even. See Table \ref{case4}. 

Notice that if $c_0$ is even, then there are no $a$-labels on the connecting points of the second-level tiles. If we observe the patterns for $c_{i-1},c_i$, then the labels of connecting edges should be $a^{(i)},b^{(i)}$. However, each $a^{(i)},b^{(i)}$ are obtained from taking sections of $a,b$. Hence the edges of all finite tiles are labeled by $a,b$. It is easy to check that on each finite tile, the labels $a,b$ are alternating. The description of the tile inflation process is complete. 

\begin{table}[H]
\centering
\caption{All cases of parities of $c_0,c_1$}
\label{allcase1}
\subcaption{Case 1: $c_0$ is odd and $c_1$ is even}
\label{case1}
\begin{adjustbox}{width=1.1\textwidth}
\begin{tabular}{|c|c|c|c|c|}
\hline
 & $\mathcal{T}_{0,1}$ & $\mathcal{T}_{1,1}$ & $\mathcal{T}_{0,2}$ & $\mathcal{T}_{1,2}$\\ \hline
Boundary points & $I_0^{(1)},I_{c_0}^{(1)}$ & $J_0^{(1)},J_{c_0-1}^{(1)}$ & $I_0^{(2)}I_{c_0}^{(1)},I_{c_0}^{(2)}J_{0}^{(1)}$ & $J_0^{(2)}I_{c_0}^{(1)},J_{C_1-1}^{(2)}J_{c_0-1}^{(1)}$  \\ \hline
\begin{tabular}{c} Connecting points\\ interchanging \\ $0,1\in V_1$ \end{tabular} & - & - & $I_0^{(2)}I_0^{(1)}\dfrac{ b}{\quad}I_1^{(2)}J_0^{(1)}$ & $J_0^{(2)}I_0^{(1)}\dfrac{ b}{\quad}J_1^{(2)}J_0^{(1)}$ \\ \hline
\begin{tabular}{c}Connecting points \\ fixing $1\in V_1$ \\ with $a$-labels \end{tabular} & - & - & \begin{tabular}{c}
     $I_{2k+1}^{(2)}J_{c_0-1}^{(1)}\dfrac{ a}{\quad}I_{2k+2}^{(2)}J_{c_0-1}^{(1)}$  \\
      for $k=0,\ldots ,c_1/2$
\end{tabular} & \begin{tabular}{c}
     $J_{2k+1}^{(2)}J_{c_0-1}^{(1)}\dfrac{a}{\quad}
J_{2k+2}^{(2)}J_{c_0-1}^{(1)}$  \\
      for $k=0,\ldots ,(c_1-2)/2$
\end{tabular} \\ \hline
\begin{tabular}{c}Connecting points \\ fixing $1\in V_1$ \\ with $b$-labels \end{tabular} & - & - & \begin{tabular}{c}
     $I_{2k}^{(2)}J_{0}^{(1)}\dfrac{ b }{\quad}
I_{2k+1}^{(2)}J_{0}^{(1)}$  \\
      for $k=1,\ldots ,c_1/2$
\end{tabular} & \begin{tabular}{c}
     $J_{2k}^{(2)}J_0^{(1)}\dfrac{ b }{\quad}J_{2k+1}^{(2)}J_0^{(1)}$  \\
      for $k=1,\ldots ,(c_1-2)/2$
\end{tabular} \\ \hline
\end{tabular}
\end{adjustbox} 

\bigskip
\subcaption{Case 2: $c_0$ is even and $c_1$ is odd}
\label{case2}
\begin{adjustbox}{width=1.1\textwidth}
\begin{tabular}{|c|c|c|c|c|}
\hline
 & $\mathcal{T}_{0,1}$ & $\mathcal{T}_{1,1}$ & $\mathcal{T}_{0,2}$ & $\mathcal{T}_{1,2}$\\ \hline
Boundary points & $I_0^{(1)},I_{c_0}^{(1)}$ & $J_0^{(1)},J_{c_0-1}^{(1)}$ & $I_0^{(2)}I_{c_0}^{(1)},I_{c_0}^{(2)}J_{c_0-1}^{(1)}$ & $J_0^{(2)}I_{c_0}^{(1)},J_{c_1-1}^{(2)}J_{c_0-1}^{(1)}$  \\ \hline
\begin{tabular}{c}Connecting points \\ interchanging \\ $0,1\in V_1$ \end{tabular}& - & - & $I_0^{(2)}I_0^{(1)}\dfrac{ b }{\quad}I_1^{(2)}J_0^{(1)}$ & $J_0^{(2)}I_0^{(1)}\dfrac{ b }{\quad}J_1^{(2)}J_0^{(1)}$ \\ \hline
\begin{tabular}{c} Connecting points \\ starting with $J_{c_0-1}^{(1)}$ \\ fixing $1\in V_1$ \end{tabular}& - & - & \begin{tabular}{c}
     $I_{2k+1}^{(2)}J_{c_0-1}^{(1)}\dfrac{ b}{\quad}I_{2k+2}^{(2)}J_{c_0-1}^{(1)}$  \\
      for $k=0,\ldots ,(c_1-1)/2$
\end{tabular} & \begin{tabular}{c}
     $J_{2k+1}^{(2)}J_{c_0-1}^{(1)}\dfrac{ b }{\quad}
J_{2k+2}^{(2)}J_{c_0-1}^{(1)}$  \\
      for $k=0,\ldots ,(c_1-3)/2$
\end{tabular} \\ \hline
\begin{tabular}{c}Connecting points \\ starting with $J_0^{(1)}$ \\ fixing $1\in V_1$ \end{tabular}& - & - & \begin{tabular}{c}
     $I_{2k}^{(2)}J_{0}^{(1)}\dfrac{ b }{\quad}
I_{2k+1}^{(2)}J_{0}^{(1)}$  \\
      for $k=1,\ldots ,(c_1-1)/2$
\end{tabular} & \begin{tabular}{c}
     $J_{2k}^{(2)}J_0^{(1)}\dfrac{ b }{\quad}J_{2k+1}^{(2)}J_0^{(1)}$  \\
      for $k=1,\ldots ,(c_1-3)/2$
\end{tabular} \\ \hline
\end{tabular}
\end{adjustbox}

\bigskip

\subcaption{Case 3: both $c_0$ and $c_1$ are odd}
\label{case3}
\begin{adjustbox}{width=1.1\textwidth}
\begin{tabular}{|c|c|c|c|c|}
\hline
 & $\mathcal{T}_{0,1}$ & $\mathcal{T}_{1,1}$ & $\mathcal{T}_{0,2}$ & $\mathcal{T}_{1,2}$\\ \hline
Boundary points & $I_0^{(1)},I_{c_0}^{(1)}$ & $J_0^{(1)},J_{c_0}^{(1)}$ & $I_0^{(2)}I_{c_0}^{(1)},I_{c_0}^{(2)}J_{c_0-1}^{(1)}$ & $J_0^{(2)}I_{c_0}^{(1)},J_{c_1-1}^{(2)}J_{0}^{(1)}$  \\ \hline
\begin{tabular}{c} Connecting points \\interchanging \\ $0,1\in V_1$ \end{tabular} & - & - & $I_0^{(2)}I_0^{(1)}\dfrac{ b }{\quad}I_1^{(2)}J_0^{(1)}$ & $J_0^{(2)}I_{0}^{(1)}\dfrac{ b }{\quad}J_1^{(2)}J_0^{(1)}$ \\ \hline
\begin{tabular}{c} Connecting points \\ fixing $1\in V_1$ \\ with $a$-labels \end{tabular} & - & - & \begin{tabular}{c}
     $I_{2k+1}^{(2)}J_{c_0-1}^{(1)}\dfrac{ a }{\quad}I_{2k+2}^{(2)}J_{c_0-1}^{(1)}$  \\
      for $k=0,\ldots ,c_1/2$
\end{tabular} & \begin{tabular}{c}
     $J_{2k+1}^{(2)}I_{c_0-1}^{(1)}\dfrac{ a }{\quad}
J_{2k+2}^{(2)}I_{c_0-1}^{(1)}$  \\
      for $k=0,\ldots ,(c_1-2)/2$
\end{tabular} \\ \hline
\begin{tabular}{c} Connecting points \\ fixing $1\in V_1$ \\ with $b$-labels \end{tabular} & - & - & \begin{tabular}{c}
     $I_{2k}^{(2)}J_{0}^{(1)}\dfrac{ b }{\quad}
I_{2k+1}^{(2)}J_{0}^{(1)}$  \\
      for $k=1,\ldots ,c_1/2$
\end{tabular} & \begin{tabular}{c}
     $J_{2k}^{(2)}J_0^{(1)}\dfrac{ b }{\quad}J_{2k+1}^{(2)}J_0^{(1)}$  \\
      for $k=1,\ldots ,(c_1-2)/2$
\end{tabular} \\ \hline
\end{tabular}
\end{adjustbox}

\bigskip
\subcaption{Case 4: both $c_0$ and $c_1$ are even}
\label{case4}

\begin{adjustbox}{width=1.1\textwidth}
\begin{tabular}{|c|c|c|c|c|}
\hline
 & $\mathcal{T}_{0,1}$ & $\mathcal{T}_{1,1}$ & $\mathcal{T}_{0,2}$ & $\mathcal{T}_{1,2}$\\ \hline
Boundary points & $I_0^{(1)},I_{c_0}^{(1)}$ & $J_0^{(1)},J_{c_0}^{(1)}$ & $I_0^{(2)}I_{c_0}^{(1)},I_{c_0}^{(2)}J_{0}^{(1)}$ & $J_0^{(2)}I_{c_0}^{(1)},J_{c_1-1}^{(2)}J_{c_0-1}^{(1)}$  \\ \hline
\begin{tabular}{c}Connecting points \\ interchanging \\ $0,1\in V_1$ \end{tabular} & - & - & $I_0^{(2)}I_0^{(1)}\dfrac{ b }{\quad}I_1^{(2)}J_0^{(1)}$ & $J_0^{(2)}I_{c_0-1}^{(1)}\dfrac{ b }{\quad}J_1^{(2)}J_0^{(1)}$ \\ \hline
\begin{tabular}{c}Connecting points \\ starting with $J_{c_0-1}^{(1)}$ \\ fixing $1\in V_1$ \end{tabular} & - & - & \begin{tabular}{c}
     $I_{2k+1}^{(2)}J_{c_0-1}^{(1)}\dfrac{ b}{\quad}I_{2k+2}^{(2)}J_{c_0-1}^{(1)}$  \\
      for $k=0,\ldots ,(c_1-1)/2$
\end{tabular} & \begin{tabular}{c}
     $J_{2k+1}^{(2)}J_{c_0-1}^{(1)}\dfrac{ b}{\quad}
J_{2k+2}^{(2)}J_{c_0-1}^{(1)}$  \\
      for $k=0,\ldots ,(c_1-3)/2$
\end{tabular} \\ \hline
\begin{tabular}{c} Connecting points \\ starting with $J_0^{(1)}$ \\ fixing $1\in V_1$ \end{tabular}& - & - & \begin{tabular}{c}
     $I_{2k}^{(2)}J_{0}^{(1)}\dfrac{ b }{\quad}
I_{2k+1}^{(2)}J_{0}^{(1)}$  \\
      for $k=1,\ldots ,(c_1-1)/2$
\end{tabular} & \begin{tabular}{c}
     $J_{2k}^{(2)}J_0^{(1)}\dfrac{\quad b \quad}{}J_{2k+1}^{(2)}J_0^{(1)}$  \\
      for $k=1,\ldots ,(c_1-3)/2$
\end{tabular} \\ \hline
\end{tabular}
\end{adjustbox}
\end{table}

By Proposition \ref{linrep} we have the following.
\begin{prop}
    The action $\Omega(\mathsf{B})\curvearrowleft D_{\infty}$ is linearly repetitive if and only if the sequence $(c_0,c_1,\ldots )$ in Theorem \ref{makebratt} is bounded. 
\end{prop}

\begin{cor}
  The dihedral group $D_{\infty}$ defined by its action on $\Omega(\mathsf{B})$ is of bounded type if and only if the sequence $(c_0,c_1,\ldots )$ is bounded. 
\end{cor}
\begin{proof}
  By the construction, there are $2$ boundary points on each finite tile. This also follows from the fact that each tile is a line segment, and thus there can be at most $2$ boundary points on each of them. There are $3$ boundary points on infinite tiles: one is fixed by $a$ and two are fixed by $b$. 
  
  By the last proposition, finite tiles are linearly repetitive. The statement is proved. 
\end{proof}

The continued fraction determined by a bounded sequence $(c_0,c_1,\ldots )$ is called a \textit{continued fraction of bounded type}. See, for instance, the discussion below Corollary 11.9 in \cite{mil2006}.

\subsection{A fragmentation of $D_{\infty}$ and its subexponential growth estimate}\label{dinftyest}
Let $(c_0,c_1,\ldots )$ be a \textbf{bounded} sequence of positive integers with each $c_i>1$. In this subsection, we describe the Grigorchuk-type fragmentation $G$ of $D_{\infty}$ determined by this sequence (with more conditions imposed) and give a subeponential growth estimate of $G$. 

By the results from the last subsection, the generator $a$ has $1$ fixed point, denoted $\mu$, and the generator $b$ has $2$ fixed points, denoted $\sigma,\lambda$. Denote by $\mu_n,\sigma_n,\lambda_n$, respectively, the $n$-th truncations of $\mu,\sigma,\lambda$. Let $n\geq 2$. Then on $\mathcal{T}_{0,n+1}$, the connector between (an isomorphic copy of) $\mathcal{T}_{0,n}$ and (an isomorphic copy of) $\mathcal{T}_{1,n}$ is of the form
\[I_{0}^{(n+1)}I_{0}^{(n)}\eta_{n-1}\dfrac{\quad x,x^{-1} \quad}{}I_1^{(n+1)}J_0^{(n)}\eta_{n-1},\addtag\label{conn1}\]
for $\eta_{n-1}\in\{\mu_{n-1},\sigma_{n-1},\lambda_{n-1}\}$ and $x\in\{a,b\}$. The edge is a two-sided arrow with labels $x$ from left to right and $x^{-1}$ backward. Similarly, in $\mathcal{T}_{1,n+1}$, the connector between (an isomorphic copy of) $\mathcal{T}_{0,n}$ and (an isomorphic copy of) $\mathcal{T}_{1,n}$ are of is form
\[J_{0}^{(n+1)}I_0^{(n)}\eta_{n-1}\dfrac{\quad x,x^{-1} \quad}{}J_1^{(n+1)}J_0^{(n)}\eta_{n-1},\addtag\label{conn2}\] 
for the same $\eta_{n-1}$ and $x$ as in (\ref{conn1}). In other words, the connectors between $\mathcal{T}_{0,n}$ and $\mathcal{T}_{1,n}$ in $\mathcal{T}_{i,n+1}$ are labeled by the same generator, for $i=0,1$.  

\begin{prop}
   Suppose $c_0$ is odd. Let $(c_{i_0},c_{i_1},c_{i_2},\ldots )$, where $i_0=0$, be the subsequence of $(c_0,c_1,\ldots )$ in which all the terms are odd. Suppose on level $n$ of $\mathsf{B}$ the connectors (\ref{conn1}),(\ref{conn2}) are labeled by $a,a^{-1}$, and the next appearance of connectors (\ref{conn1}),(\ref{conn2}) labeled by $a,a^{-1}$ is on level $n+l_{a,n,\mu}$. Then $l_{a,n,\mu}$ is uniformly bounded if and only if $|i_k-i_{k-1}|$ is uniformly bounded, for $k\geq 0$. 
\end{prop}
\begin{proof}
   On level $1$, these $2$ connectors are always labeled by $a$. From level $1$ to $i_1$, the generator $a$ is directed along and fixes the path $\mu_{i_1}=J_{c_{i_1-1}-1}^{(i_1)}\ldots J_{c_1-1}^{(2)}J_{c_0-1}^{(1)}$. The section of $a$ on $\gamma$ is $b^{(i_1+1)}$ which has nontrivial section on $I_{c_{i_1}}^{(i_1+1)}$. It follows that on level $c_{i_1+2}$ the connectors (\ref{conn1}),(\ref{conn2}) are labeled by $a,a^{-1}$. The argument for the levels from $i_1$ to $i_2$ follows similarly. The proof is complete. 
\end{proof}

By the proposition above, we add the following assumption to the sequence $(c_0,c_1,\ldots )$. 
\begin{Assumptions}[leftmargin=5\parindent]
    \item[\textbf{Assumption:}] The consecutive appearances of odd terms of $(c_0,c_1,\ldots )$ have uniformly bounded gaps. 
\end{Assumptions}

It is also observed directly from the construction that the consecutive appearances of the connectors (\ref{conn1}),(\ref{conn2}) labeled by $b,b^{-1}$ with $\eta_{n-1}=\sigma_{n-1}$ have uniformly bounded gaps. Denote this gap by $l_{b,n,\sigma}$. The same holds for labels $b,b^{-1}$ with  $\eta_{n-1}=\lambda_{n-1}$. Denote this gap by $l_{b,n,\lambda}$. 

The support of $a$ denoted Supp$(a)=\Omega(\mathsf{B})\backslash\{\mu\}$, and Supp$(b)=\Omega(\mathsf{B})\backslash\{\sigma,\lambda\}$. Let us partition subsets of Supp$(a)$ and Supp$(b)$ such that they accumulate respectively on $\mu,\sigma,\lambda$. For convenience, let $\mathsf{B}_{a,\mu}$ be the Bratteli diagram obtained from a telescoping of bounded steps of $\mathsf{B}$ such that the consecutive connectors of the forms (\ref{conn1}),(\ref{conn2}) are labeled by $a,a^{-1}$. Similarly define $\mathsf{B}_{b,\sigma},\mathsf{B}_{b,\lambda}$ to be, respectively, telescopings of bounded steps of $\mathsf{B}$ such that the consecutive appearances of connectors (\ref{conn1}),(\ref{conn2}) have labels $b,b^{-1}$ and $\eta_{n-1}$ are respectively $\sigma_{n-1}$ and $\lambda_{n-1}$. The telescoped $\mu,\sigma,\lambda$, their truncations, and the connecting points in the connectors (\ref{conn1}),(\ref{conn2}) are denoted the same as above. Let $\mathsf{s}$ be the shift map as before. Define 
\begin{align*}
      W_{x,\eta,n}= (\mathsf{s}^{n+1}(\Omega(\mathsf{B}_{x,\eta}))I_{0}^{(n+1)}I_0^{(n)}\eta_{n-1})\cup(\mathsf{s}^{n+1}(\Omega(\mathsf{B}_{x,\eta}))I_{0}^{(n+1)}J_0^{(n)}\eta_{n-1})\\
      \cup(\mathsf{s}^{n+1}(\Omega(\mathsf{B}_{x,\eta}))J_{0}^{(n+1)}I_0^{(n)}\eta_{n-1})\cup(\mathsf{s}^{n+1}(\Omega(\mathsf{B}_{x,\eta}))J_{1}^{(n+1)}J_0^{(n)}\eta_{n-1}),
\end{align*}
for $x\in\{a,b\}$ and $\eta\in\{\mu,\sigma,\lambda\}$. Then $W_{a,\mu}$ is $a$-invariant and $W_{b,\sigma},W_{b,\lambda}$ are $b$-invariant. Let $P_{\mu,i}=\bigcup_{k=0}^{\infty} W_{\mu,3k+i}$, for $i=0,1,2$. Each $P_{\mu,i}$ accumulates on $\mu$. Then the collection $\mathcal{P}_{\mu}=\{P_{\mu,0},P_{\mu,1},P_{\mu,2}\}$ forms an open partition of the subset containing only the connectors (\ref{conn1}),(\ref{conn2}) of Supp$(a)$. Denote this subset by $A_{\mu}$. Let $b_0,c_0,d_0$ be the homeomorphisms acting trivially on $P_{\mu,2},P_{\mu,1},P_{\mu,0}$, respectively, and as $a$ on their complements in $A_{\mu}$. Similarly, we can define, for $i=0,1,2$, $P_{\sigma,i}=\bigcup_{k=1}^{\infty}W_{\sigma,3k+i}$ and $P_{\lambda,i}=\bigcup_{k=1}^{\infty}W_{\lambda,3k+i}$, and the subsets of Supp$(b)$, denoted $B_{\sigma},B_{\lambda}$, such that the collection $\mathcal{P}_{\sigma}$ is an open partition of $B_{\sigma}$ accumulating on $\sigma$ and $\mathcal{P}_{\lambda}$ is an open partition of $B_{\lambda}$ accumulating on $\lambda$. Let $b_1,c_1,d_1$ act as identity respectively on $P_{\sigma,2},P_{\sigma,1},P_{\sigma,0}$ and as $a$ on their complements in $B_{\sigma}$; let $b_2,c_2,d_2$ act as identity respectively on $P_{\lambda,2},P_{\lambda,1},P_{\lambda,0}$ and as $a$ on their complements in $B_{\lambda}$. Let $a',b'$ act as $a,b$ on the complement of $A_{\mu}\cup B_{\sigma}\cup B_{\lambda}$. The fragmented group $G$ is then generated by $a',b',b_i,c_i,d_i$, for $i=0,1,2$. It follows that the points $\mu,\sigma,\lambda$ are purely non-Hausdorff singularities of $G$. Relabel them respectively by $\xi_0,\xi_1,\xi_2$. Their respective orbital graphs $\Gamma_{\xi_i}$ are, for $i=0, 1, 2$, the inductive limits of embeddings of $\mathcal{T}_{1,3k+i}$ to the end of $\mathcal{T}_{1,3(k+1)+i}$ that contains $\xi_{i,3(k+1)+i}$ (the $(3(k+1)+i)$-th truncation of $\xi_i$), and eventually adding loop at the boundary point of $\mathcal{T}_{\xi_i}$. The graph of germs $\widetilde{\Gamma}_{\xi_i}$ is obtained by connecting four copies of $\Gamma_{\xi_i}$ by the Cayley
graph of the four-group $H_i=\{Id,b_i,c_i,d_i\}$ at $\xi_i$. 

 Let $g\in G$. A \textit{rank $n$ traverse of type $h$} is a walk of $g$ starting and ending at a boundary points of $\mathcal{T}_{j,n}$, not touching the boundary points in between, and this walk can be lifted to a central part of a $\widetilde{\Gamma}_{\xi_i}$ starting in the branch $(Id,\Gamma_{\xi_i})$ and ending in the branch $(h,\Gamma_{\xi_i})$. Denote by $\Theta_{n,h}$ the set of rank $n$ traverse of type $h$ and let $\Theta_n=\bigcup\limits_{h\in H_{i}}\Theta_{n,h}$. Let $T_{n,h}=\#\Theta_{n,h}$, and $T_n=\sum\limits_{h\in H_{i}\backslash\{Id\}}T_{n,h}$. Let $F_h(t)=\sum\limits_{n=0}^{\infty}T_{n,h}t^n$ and $F(t)=\sum\limits_{n=0}^{\infty}T_nt^n$. Since the gaps $l_{a,n,\mu},l_{b,n,\sigma},l_{b,n,\lambda}$ are uniformly bounded, there exists $N\in\mathbb{N}$ such that they are bounded above by $N$. It follows that the tiles on level at least $n+(3N+3)$ contains all types of connectors of tiles on levels $n$. Let $L=3N+3$. Then a traverse of $\mathcal{T}_{j,n+L}$ must induce all types of traverses of all $\mathcal{T}_{j,n}$. It follows that 
 \[T_{n+L}\leq \sum_{h\neq Id}T_{n,h},\]
\[\implies t^{-L}(F(t)-(T_0+T_1t+T_2t^2+\ldots +T_{L-1}t^{L-1}))\leq \sum\limits_{h\neq Id}F_h(t)\]
\[\implies F(t)\leq |w| \dfrac{g(t)}{Ct^{-L}-D}\]
where $g(t)>0$ for $t>0$ and $0<D<C$. Hence if $t\in (0,\eta^{-1})$, where $\eta=\sqrt[L]{\frac{D}{C}}$ is the real positive root of the polynomial equation $Cx^L-D=0$, the right-hand side of the above inequality is positive.

Now let us estimate the cardinalities $|\mathcal{T}_{j,n}|$. Since the sequence $(c_0,c_1,\ldots )$ is bounded, there exists $M>1$ such that each $c_i\leq M$. Then we have
\[|\mathcal{T}_{0,n+1}|\leq |\mathcal{T}_{0,n}|+M|\mathcal{T}_{1,n}|\] and \[|\mathcal{T}_{1,n+1}|\leq |\mathcal{T}_{0,n}|+(M-1)|\mathcal{T}_{1,n}|<|\mathcal{T}_{0,n}|+M|\mathcal{T}_{1,n}|.\]
Hence $\limsup_{n\rightarrow\infty}(\sum_{j=0}^1|\mathcal{T}_{j,n}|)^{1/n}$ exists and is denoted $\beta$. Note that $\beta>1$. It follows that the growth function of $G$, denoted $\gamma_G(R)$ is dominated by $\exp(R^{\alpha})$ for every $\alpha>\frac{\log \beta}{\log \beta-\log\eta}$. 

\section{A family of groups with non-linear orbital graphs}\label{nonlin13}
\subsection{The construction}
Let $w=\ldots w_n\ldots w_2w_1$ be a left-infinite sequence of $1,3$, i.e., $w\in\{1,3\}^{\omega}$. Consider the Bratteli diagram $\mathsf{B}_w$ determined by $w$: if $w_n=1$, then level $n$ of $\mathsf{B}_w$ is shown in Figure \ref{llv1}, and if $w_n=3$, then level $n$ of $\mathsf{B}_w$ is shown in Figure \ref{llv3}. For simplicity, we also impose the following conditions on $w$:
\begin{enumerate}
\item the subword $33$ is prohibited in $w$;
\item\label{bdd33} the consecutive appearances of $3$ have bounded gaps.  

\end{enumerate}

\begin{figure}
    \centering

    % First subfigure (left)
    \begin{subfigure}{0.45\textwidth}
        \centering
        \begin{adjustbox}{width=.9\textwidth}
            \begin{tikzpicture}
                \begin{scope}state without output/.append style={draw=none}%[every node/.style={circle,thick,draw}]
                    \node (A) at (0,0) {};
                    \node (B) at (3,0) {};
                    \node (C) at (0,-3) {};
                    \node (D) at (3,-3) {};
                \end{scope}

                \begin{scope}
                    \path [-] (A) edge node [midway,left] {0} (C);
                    \path [-] (A) edge node [near end,right]{$e_1$} (D);
                    \path [-] (B) edge node [near start,left] {$e_2$} (C);

                    \filldraw[black] (0,0) circle  (2.5pt);
                    \filldraw[black] (3,0) circle (2.5pt);
                    \filldraw[black] (0,-3) circle (2.5pt);
                    \filldraw[black] (3,-3) circle (2.5pt);
                \end{scope}
            \end{tikzpicture}
        \end{adjustbox}
        
        \caption{Level corresponding to $w_n=1$.}
        \label{llv1}
    \end{subfigure}
    \hfill
    % Second subfigure (right)
    \begin{subfigure}{0.45\textwidth}
        \centering
        \begin{adjustbox}{width=1\textwidth}
            \begin{tikzpicture}
                \begin{scope}state without output/.append style={draw=none}%[every node/.style={circle,thick,draw}]
                    \node (A) at (0,0) {};
                    \node (B) at (3,0) {};
                    \node (C) at (0,-3) {};
                    \node (D) at (3,-3) {};
                \end{scope}

                \begin{scope}
                    \path [-] (A) edge [bend right=20] node [midway,left] {$0$} (C);
                    \path [-] (A) edge node [midway,right] {$1$} (C);
                    \path [-] (A) edge [bend left=20] node [midway,right] {$2$} (C);
                    \path [-] (B) edge node [near start,left] {$e_2$} (C);
                    \path [-] (A) edge node [near end,right]{$e_1$} (D);

                    \filldraw[black] (0,0) circle  (2.5pt) node[anchor=west] {};
                    \filldraw[black] (3,0) circle (2.5pt) node[anchor=west] {};
                    \filldraw[black] (0,-3) circle (2.5pt) node[anchor=west] {};
                    \filldraw[black] (3,-3) circle (2.5pt) node[anchor=west] {};
                \end{scope}
            \end{tikzpicture}
        \end{adjustbox}

        \caption{Level corresponding to $w_n=3$.}
        \label{llv3}
    \end{subfigure}

    \caption{Bratteli Diagram $\mathsf{B}_w$}
    \label{exmpnonlin}
\end{figure}

Let $S=\{a,b,c\}$ where $a,b$ are involutions and $c$ has order $3$. Let $G_0$ be the group generated by $S$ defined by the following tile inflation rules. Denote by $\mathcal{T}_{1,n}$ and $\mathcal{T}_{2,n}$ the tiles on level $n$ of $\mathsf{B}_w$ corresponding to the vertex on the left and on the right, respectively. Tiles on levels $1,2$ are shown below, where only active edges are displayed. 

\refstepcounter{casecounter} % Increment the counter
\textbf{Case \thecasecounter}:\label{case01} $w_1=w_2=1$, then

\begin{align*}
    &\mathcal{T}_{1,1}:\text{  }0\frac{ \ \quad b \quad \ }{}e_2,\text{       }\mathcal{T}_{1,2}:\text{  }e_1,\\
     \mathcal{T}_{2,1}:\text{  }e_2e_1&\frac{ \ \quad a \quad \ }{}00\frac{ \ \quad b \quad \ }{}0e_2,\text{       }\mathcal{T}_{2,2}:\text{  }e_10\frac{ \ \quad b \quad \ }{}e_1e_2=e_1\mathcal{T}_{1,1}.
\end{align*}
The edges labeled by $a,b$ above are active, and they are two-sided arrows. 

\refstepcounter{casecounter}
\textbf{Case \thecasecounter}:\label{case02} $w_1=3$ and $w_2=1$, then 
\begin{align*}
    &\mathcal{T}_{1,1}:\text{ see Figure \ref{t1131} },\text{       }\mathcal{T}_{1,2}:\text{  }e_1,\\
     &\mathcal{T}_{2,1}:\text{ see Figure \ref{t2131} },\text{       }\mathcal{T}_{2,2}:\text{  }e_1\mathcal{T}_{1,1}.
\end{align*}

\begin{figure}[htb]
\centering
 \begin{subfigure}{0.45\textwidth}
 \centering
\begin{tikzpicture}
      % Nodes
    \node[fill, circle, inner sep=2pt, label=below:{$2$}] (A) at (0,0) {};
    \node[fill, circle, inner sep=2pt, label=left:{$0$}] (B) at (-2,1) {};
    \node[fill, circle, inner sep=2pt, label=left:{$1$}] (D) at (-2,-1) {};
    \node[fill, circle, inner sep=2pt, label=above:{$e_2$}] (F) at (1.5,0) {};

    % Edges
    \draw (A) -- node[above] {$b$} (F);
    \draw [->] (A) -- node[left] {$c$} (B);
    \draw (B) -- node[left] {$c$} (A);
    \draw [->] (D) -- node[right] {$c$} (A);
    \draw [->] (B) -- node[left] {$c$} (D);
    
\end{tikzpicture}
\caption{$\mathcal{T}_{1,1}$.}
\label{t1131}
\end{subfigure}
\vspace{.5cm}
 \begin{subfigure}{0.45\textwidth}
\centering
\begin{tikzpicture}
      % Nodes
    \node[fill, circle, inner sep=2pt, label=below:{$02$}] (A) at (0,0) {};
    \node[fill, circle, inner sep=2pt, label=left:{$00$}] (B) at (-2,1) {};
    \node[fill, circle, inner sep=2pt, label=left:{$01$}] (D) at (-2,-1) {};
    \node[fill, circle, inner sep=2pt, label=above:{$0e_2$}] (F) at (1.5,0) {};
     \node[fill, circle, inner sep=2pt, label=left:{$e_2e_1$}] (E) at (-3,-2) {};

    % Edges
    \draw (A) -- node[above] {$b$} (F);
    \draw [->] (A) -- node[left] {$c$} (B);
    \draw (B) -- node[left] {$c$} (A);
    \draw [->] (D) -- node[right] {$c$} (A);
    \draw [->] (B) -- node[left] {$c$} (D);
    \draw (D) -- node[below] {$a$} (E);
    
\end{tikzpicture}
\caption{$\mathcal{T}_{2,1}.$}
\label{t2131}
\end{subfigure}
\caption{Tiles when $w_1=3$ and $w_2=1$.}
\end{figure}

\refstepcounter{casecounter}
\textbf{Case \thecasecounter}:\label{case03} $w_1=1$ and $w_2=3$, then

\begin{align*}
    &\mathcal{T}_{1,1}:\text{  }0\frac{ \ \quad b \quad \ }{}e_2,\text{       }\mathcal{T}_{1,2}:\text{  }e_1,\\
     \mathcal{T}_{2,1}: &\text{    see Figure \ref{t2113} }, \text{       }\mathcal{T}_{2,2}:\text{  }e_10\frac{ \ \quad b \quad \ }{}e_1e_2e_1=\mathcal{T}_{1,1}.
\end{align*}  

\begin{figure}
\centering
\begin{tikzpicture}
      % Nodes
    \node[fill, circle, inner sep=2pt, label=below:{$20$}] (A) at (0,0) {};
    \node[fill, circle, inner sep=2pt, label=left:{$00$}] (B) at (-2,1) {};
    \node[fill, circle, inner sep=2pt, label=left:{$10$}] (D) at (-2,-1) {};
    \node[fill, circle, inner sep=2pt, label=left:{$0e_2$}] (C) at (-3,2.5) {};
    \node[fill, circle, inner sep=2pt, label=left:{$1e_2$}] (E) at (-3,-2.5) {};
    \node[fill, circle, inner sep=2pt, label=above:{$2e_2$}] (F) at (1.5,1.5) {};
    \node[fill, circle, inner sep=2pt, label=right:{$e_2e_1$}] (G) at (4,0) {};

    % Edges
    \draw (A) -- node[above] {$a$} (G);
    \draw (A) -- node[above] {$b$} (F);
    \draw [->] (A) -- node[left] {$c$} (B);
    \draw (B) -- node[left] {$b$} (C);
    \draw (B) -- node[left] {$c$} (A);
    \draw (D) -- node[left] {$b$} (E);
    \draw [->] (D) -- node[right] {$c$} (A);
    \draw [->] (B) -- node[left] {$c$} (D);
    
\end{tikzpicture}
\caption{$\mathcal{T}_{2,1}$ when $w_1=1$ and $w_2=3$.}
\label{t2113}
\end{figure}

The boundary points on infinite tiles of $a,b,c$ are described as follows. Let $s_{n}\in\{0,1\}$ and $t_{n}\in\{0,2\}$. The generator $a$ has $2$ boundary points $\ldots e_2e_1s_{3n-2}\ldots e_2e_1s_4e_2e_1s_1$ and $\ldots e_1s_{3n-1}e_2\ldots e_1s_5e_2e_1s_2e_2$, which can be formally written respectively as $\sigma=(e_2e_1s_{3n-2})^{\omega}$ and $\lambda=(e_2s_{3n-1}e_1)^{\omega}$ for $n\in\mathbb{N}$, noting that they are not periodic sequences. The letter $s_k=0$ if $w_k=1$ and $s_k=1$ if $w_k=3$. The generator $b$ has $1$ boundary point, written formally as $\mu=(t_{3n}e_2e_1)^{\omega}$ for $n\in\mathbb{N}$. Similarly, the letter $t_k=0$ if $w_k=1$ and $t_k=2$ if $w_k=3$. The generator $c$ has $1$ boundary point, denoted $\delta=\ldots x_2x_1$, where each $x_i$ is described by the following rule (note that $33$ is prohibited in $w$):
\begin{align*}
 x_i=0 &\text{ if } w_i=1,\\
 x_i=e_1, \text{ } x_{i+1}=e_2 &\text{ if } w_i=3.
\end{align*}
All the generators fix their boundary points on infinite tiles. 

Now let us define the connectors of tiles on each level of $\mathsf{B}_{w}$. We use triple notations (see Subsection \ref{tia}) to represent the connectors labeled by the involutions $a,b$ (as two-sided arrows), and ($3$-)cycle notations $(\alpha,\beta,\gamma)$ to represent the connectors labeled by $c$, where $\alpha\cdot c=\beta$,  $\beta\cdot c=\gamma$, and $\gamma\cdot c=\alpha$. Since each $\mathcal{T}_{2,n+1}$ is an isomorphic copy of $\mathcal{T}_{1,n}$, we only need to define the connectors on each $\mathcal{T}_{1,n}$. 

For $n=1$, the connectors on $\mathcal{T}_{1,1}$ are $(0,e_2,b)$ in Cases \ref{case01}, \ref{case03}, and $(0,e_2,b)$, $(0,1,2)$ in Case \ref{case02}. For $n=2$, the connectors are $(e_2e_1,00,a)$ in Case \ref{case01}, $(01,e_2e_1,a)$ in Case \ref{case02}, and $(20,e_1e_2,a)$, $(00,10,20)$ in Case \ref{case03}. 

For $n\geq 3$, $c$ is only active when $w_n=3$, and the $c$-connectors are of the forms
\begin{align}
    (000\delta_{n-3},100\delta_{n-3},200\delta_{n-3}) \text{ if } w_{n-2}=1,\label{order3con1}\\
  (0e_2e_1\delta_{n-3},1e_2e_1\delta_{n-3},2e_2e_1\delta_{n-3}) \text{ if } w_{n-2}=3. \label{order3con2}
\end{align}
If $n=3$, then $\delta_{n-3}$ is the empty word. 

The connectors $a,b$ are defined as follows. Let $x_n\in\{a,b\}$ be the labels. Then, 
\begin{align}\label{connectorab}
x_n=
    \begin{cases}
    b\text{, if } n\equiv1 \pmod{3},\\
    a\text{, otherwise. }
    \end{cases}
\end{align}
The tile $\mathcal{T}_{1,n}$ is partly obtained by taking $1$ isomorphic copy $y\mathcal{T}_{1,n-1}$ ($y=s_n$ or $t_n$) of $\mathcal{T}_{1,n-1}$ and $1$ isomorphic copy $e_2\mathcal{T}_{2,n-1}$ of $\mathcal{T}_{1,n-1}$, and connecting them at the connecting points by $x_n$. The connectors are of the form $(t_{3k+1}t_{3k}\mu_{3k-1},e_2e_1\mu_{3k-1},b)$, $(s_{3k+2}s_{3k+1}\sigma_{3k},e_2e_1\sigma_{3k},a)$, and $(s_{3k+3}s_{3k+2}e_2\lambda_{3k},e_2e_1e_2\lambda_{3k},a)$. See Table \ref{fibobdry} for the locations of $a,b$-connectors and continuations of boundary points of $a,b$ on finite tiles. 

Note that we only defined $a,b,c$ along their boundary points. Extending them identically to the rest of $\Omega(\mathsf{B}_w)$, the inverse semigroup $G_0$ generated by $\{a,b,c\}$ is a group. The construction is complete. 

\begin{table}
\centering
\begin{tabular}{|c|c@{\hskip 1in}c|}
\hline
    Level $n$ & $\mathcal{T}_{1,n}$ & $\mathcal{T}_{2,n}$ \\ \hline\hline
    $3k-1$ & $\tikzmark{a1}\mu_{3k-1}$\tikzmark{a} & \\ [1cm]
    $3k$   & \tikzmark{b1}$\mu_{3k}$ \tikzmark{b}  & \tikzmark{c}$e_1\mu_{3k-1}$ \\[1cm]
    $3k+1$ & \tikzmark{d1}$t_{3k+1}t_{3k}\mu_{3k-1}\dfrac{b}{\quad}e_2e_1\mu_{3k-1}$\tikzmark{d} & $\tikzmark{e}\mu_{3k+1}$\\[1cm]
    $3k+2$ & $\mu_{3k+2}$ \tikzmark{f}&\\
\hline \hline  
    $3k-1$ &              & \tikzmark{g}$\sigma_{3k-1}$ \\ [1cm]
    $3k$   & \tikzmark{h1}$\sigma_{3k}$\tikzmark{h} &  \\[1cm]
    $3k+1$ & \tikzmark{i1}$s_{3k+1}\sigma_{3k}=\sigma_{3k+1}$\tikzmark{i} & \tikzmark{j}$e_1\sigma_{3k}$ \\[1cm]
    $3k+2$ & \tikzmark{k1}$s_{3k+2}s_{3k+1}\sigma_{3k}\dfrac{a}{\quad}e_2e_1\sigma_{3k}$\tikzmark{k} & \tikzmark{l}$\sigma_{3k+2}$\\
\hline\hline
    $3k-1$ & $\lambda_{3k-1}$\tikzmark{m} & \\ [1cm]
    $3k$   &  & $\tikzmark{n}\lambda_{3k}$ \\[1cm]
    $3k+1$ & \tikzmark{o1}$e_2\lambda_{3k}=\lambda_{3k+1}$\tikzmark{o} & \\[1cm]
    $3k+2$ & \tikzmark{p1}$s_{3k+2}e_2\lambda_{3k}=\lambda_{3k+2}$\tikzmark{p} & \tikzmark{q}$e_1e_2\lambda_{3k}$\\[1cm]
    $3k+3$ & \tikzmark{r1}$s_{3k+3}s_{3k+2}e_2\lambda_{3k}\dfrac{a}{\quad} e_2e_1e_2\lambda_{3k}$\tikzmark{r} & \tikzmark{s}$\lambda_{3k+3}$\\
\hline    
  \end{tabular}
  \begin{tikzpicture}[overlay, remember picture, shorten >=.5pt, shorten <=.5pt, transform canvas={yshift=.25\baselineskip}]
  \draw [->] ({pic cs:a1}) [bend right] to node [left] {$t_{3k}$} ({pic cs:b1}) ;
    \draw [->] ({pic cs:a}) to node [right] {$e_1$} ({pic cs:c});
    \draw [->] ({pic cs:b1}) [bend right] to node [right] {$t_{3k+1}$} ({pic cs:d1});
    \draw [->] ({pic cs:c}) to node [near start, right] {$e_2$} ({pic cs:d});
    \draw [->] ({pic cs:b}) to node [near start, right] {$e_1$} ({pic cs:e});
    \draw [->] ({pic cs:e}) to node [right] {$e_2$} ({pic cs:f});

    \draw [->] ({pic cs:g}) to node [left] {$e_2$} ({pic cs:h});
    \draw [->] ({pic cs:h1}) [bend right] to node [right] {$s_{3k+1}$} ({pic cs:i1});
    \draw [->] ({pic cs:i1}) [bend right] to node [right] {$s_{3k+2}$} ({pic cs:k1});
    \draw [->] ({pic cs:h}) to node [right] {$e_1$} ({pic cs:j});
    \draw [->] ({pic cs:i}) to node [near start, right] {$e_1$} ({pic cs:l});
    \draw [->] ({pic cs:j}) to node [near start, right] {$e_2$} ({pic cs:k});

    \draw [->] ({pic cs:m}) to node [right] {$e_1$} ({pic cs:n});
    \draw [->] ({pic cs:n}) to node [right] {$e_2$} ({pic cs:o});
    \draw [->] ({pic cs:o1}) [bend right] to node [right] {$s_{3k+2}$} ({pic cs:p1});
    \draw [->] ({pic cs:p1}) [bend right] to node [right] {$s_{3k+3}$} ({pic cs:r1});
    \draw [->] ({pic cs:o}) to node [right] {$e_1$} ({pic cs:q});
    \draw [->] ({pic cs:q}) to node [near start, right] {$e_2$} ({pic cs:r});
    \draw [->] ({pic cs:p}) to node [near start, right] {$e_1$} ({pic cs:s});
  \end{tikzpicture}
  
  \caption{Locations and continuations of boundary points and locations of connecting points of $a,b$, for $k\geq 1$.}
  \label{fibobdry}
\end{table}

\subsection{A fragmentation of $G_0$}
We fragment $\langle a \rangle$, $\langle b \rangle$ together and $\langle c\rangle$ independently. The first fragmentation is similar to that of the Golden Mean dihedral group described in \cite[Subsection 4.2]{BNZ}, and the second fragmentation is similar to that of the Fabrykowski-Gupta group described in \cite[Subsection 5.4]{kua24a}. 

\textbf{Now let us fragment $\langle a \rangle$, $\langle b \rangle$.} Define $W_{\mu,n}$, for $n\geq 1$, to be the union of the sets of paths in $\Omega(\mathsf{B}_w)$ starting with $t_{3n+1}t_{3n}e_2e_1(t_{3k}e_2e_1)^{n-1}$ and $e_2e_1e_2e_1(t_{3k}e_2e_1)^{n-1}$, $k=1,\ldots ,n-1$. Define $W_{\mu,0}$ to be the union of the sets of paths starting with $t_1$ and $e_2$. Then $W_{\mu,n}$, for $n\geq 0$, are $b$-invariant. Let $P_{\mu,i}=\bigcup_{k=0}^\infty W_{\mu,3k+i}$ for $i=0,1,2$. The collection $\mathcal{P}_{\mu}=\{P_{\mu,0},P_{\mu,1},P_{\mu,2}\}$ is an open partition of $\textup{Supp}(b)\backslash \{\mu\}$. Let $b_0,c_0,d_0$ be the homeomorphisms acting trivially on $P_2,P_1,P_0$, respectively, and as $b$ on their complements. Then we have the fragmentation $F_b=\langle b_0,c_0,d_0 \rangle$. Similarly, we can fragment $a$ around its fixed points $\sigma=(e_2e_1s_k)^{\omega}$ and $\lambda=(e_1s_{k+1}e_2)^{\omega}$. Define $W_{\sigma,n}$ to be the union of the sets of paths in $\Omega(\mathsf{B}_w)$ starting with $s_{3n+2}s_{3n+1}(e_2e_1s_k)^{n}$ and $e_2e_1(e_2e_1s_k)^{n}$; define $W_{\lambda,n}$ to be the union of the sets of paths in $\Omega(\mathsf{B})$ starting with $s_{3n+3}s_{3n+2}e_2(e_1s_{k+1}e_2)^{n}$ and $e_2e_1e_2(e_1s_{k+1}e_2)^{n}$, for $n\geq 0$ and $k=1,\ldots ,n$. Both $W_{\sigma,n}$ and $W_{\lambda,n}$ are $a$-invariant. For $i=0,1,2$, let $P_{\sigma,i}=\bigcup_{k=1}^{\infty}W_{\sigma,3k+i}$ and $P_{\lambda,i}=\bigcup_{k=1}^{\infty}W_{\lambda,3k+i}$. The collection $\mathcal{P}_{\sigma,\lambda}=\{P_{\sigma,i},P_{\lambda,i}\}_{i=0}^2$ forms an open partition of $\textup{Supp}(a)\backslash\{\sigma,\lambda\}$. Let $b_1,c_1,d_1$ act as identity respectively on $P_{\sigma,2},P_{\sigma,1},P_{\sigma,0}$ and as $a$ on their complements; let $b_2,c_2,d_2$ act as identity respectively on $P_{\lambda,2},P_{\lambda,1},P_{\lambda,0}$ and as $a$ on their complements. Moreover, define $a_0,a_1,a_2$ as follows. For $\gamma\in\Omega(\mathsf{B}_w)$, 
\begin{align*}
    &\gamma s_1\cdot a_0=e_2, \text{  }\gamma e_2\cdot a_0=\gamma s_1, \text{  } \gamma e_1\cdot a_0=\gamma e_1,\\
    \gamma s_2s_1  \cdot &a_1=\gamma e_2e_1,\text{  } \gamma e_2e_1\cdot a_1=\gamma s_2s_1,\text{  } \gamma e_2\cdot a_1=\gamma e_2,\\
    \gamma s_4s_3e_2  \cdot &a_2=\gamma \cdot e_2e_1e_2,\text{  } \gamma e_2e_1e_2\cdot a_2=\gamma s_4s_3e_2,\text{  } \gamma e_1\cdot a_2=\gamma e_1.
\end{align*}
Let $G_{a,b}$ be the group generated by $a_i,b_i,c_i,d_i$, for $i=0,1,2$. The generators are extended identically to the rest of $\Omega(\mathsf{B}_w)$ if they are not defined on the whole $\Omega(\mathsf{B}_w)$.  This completes the fragmentation. 

\textbf{Now let us fragment $\langle c \rangle$.} To simplify the description, we perform a telescoping of bounded steps of $\mathsf{B}_w$ such that the new diagram is corresponding to the word $w=\ldots 131313$. This is possible due to Condition \ref{bdd33}. Denote this new Bratteli diagram by $\mathsf{B}_c$. Then the boundary point on infinite tile of $c$ is $\delta=(e_2e_1)^{\omega}\in\Omega(\mathsf{B}_c)$. 

Define $W_{2n}$ to be the set of infinite paths in $\mathsf{B}_c$ starting with $0\delta_{2n}$, $1\delta_{2n}$, $2\delta_{2n}$. Then Supp$(c)=\bigcup\limits_{n=0}^{\infty}W_{2n}$.  Let $P_0=\bigcup\limits_{k=0}^{\infty}W_{2\cdot 4k}$, $P_1=\bigcup\limits_{k=0}^{\infty}W_{2\cdot(4k+1)}$, $P_2=\bigcup\limits_{k=0}^{\infty}W_{2\cdot(4k+2)}$, and $P_3=\bigcup\limits_{k=0}^{\infty}W_{2\cdot(4k+3)}$. Then $\mathcal{P}=\{P_0,P_1,P_2,P_3\}$ forms a partition of Supp$(c)$, and each piece is $c$-invariant accumulating on $\delta$. Let $f_1$ act as $c$ on $P_1\cup P_3$, as $c^2$ on $P_2$, and as identity on $P_0$. Let $f_2$ act as $c$ on $P_0\cup P_1\cup P_2$ and as identity on $P_3$. Define $G_c=\langle f_1,f_2 \rangle$. This fragmentation can also be described by the \textit{subdirect product}\footnote{A subdirect product is an injective homomorphism (embedding) $\pi=(\pi_1,...,\pi_n):H\rightarrow (\mathbb{Z}/k\mathbb{Z})^n$ in which each $\pi_i$ is surjective, for $i\in\{1,...,n\}$. See \cite[Chapter 4]{JC19} for a detailed discussion and \cite[Subsection 2.4]{kua24a} for a summary.} $G_c\rightarrow \langle c \rangle^4\cong (\mathbb{Z}/3\mathbb{Z})^4$ defined by 
 \begin{align*}
    \begin{cases}
    f_1\mapsto (0,1,2,1),\\
    f_1\mapsto (1,1,1,0).\\
    \end{cases}
\end{align*}
In other words, we have epimorphisms $\pi_i:G_c\rightarrow \mathbb{Z}/3\mathbb{Z}$ defined by $P_i$, for $i=0,1,2,3$. The connectors are shown in Figure \ref{figconc}, and they appear in the order of (1)(2)(3)(4) along the truncations of $\delta$ at the (telescoped) connecting points described in (\ref{order3con1}) and (\ref{order3con2}). In other words, the connecting points are of the form $0\delta_{2n},1\delta_{2n},2\delta_{2n}$, for $n\geq 0$. The fragmentation is complete. 

\begin{figure}
\centering
\begin{minipage}{.5\textwidth}
\begin{tikzpicture}
\begin{scope}state without output/.append style={draw=none}%[every node/.style={circle,thick,draw}]
    \node (A) at (3,0) {};
    \node (B) at (1,-3) {};

    \node (C) at (5,-3) {};
   
\end{scope}

\begin{scope} %[>={Stealth[black]},
              every node/.style={fill=white,circle},
              every edge/.style={draw=black,very thick}]
              
    %\path [->] (A) edge [bend left=40] node {} (B);
    %\path [->] (B) edge [bend left=40] node {} (A);
    %\draw [->] (B) to  [out=320,in=40,looseness=9] (B);
    %\path [->] (D) edge node {$3$} (C);
    %\path [->] (A) edge node {$3$} (E);
    %\path [->] (D) edge node {$3$} (E);
    %\path [->] (D) edge node {$3$} (F);
    %\path [->] (C) edge node {$5$} (F);
    %\path [->] (E) edge node {$8$} (F); 
    %\path [->] (B) edge[bend right=60] node {$1$} (E); 

    %\draw[]  (A) node[draw=none][midway,above] {$b_1$} (B);
    \path [->] (A) edge node [midway,above] {$f_2$} (B);
   
    \path [->] (B) edge node [midway,above] {$f_2$} (C);
    \path [->] (C) edge node [midway,left] {$f_2$} (A);
   
    %\draw [->] (A) to  [out=90,in=140,looseness=8] edge node [midway,left] {$c_2$} (A);
    \path [->] (A) edge [out=30, in=150, looseness=15, "$f_1$"] (A);
     \path [->] (B) edge [out=150, in=250, looseness=15, "$f_1$"] (B);
      \path [->] (C) edge [out=310, in=40, looseness=15, "$f_1$"] (C);
      
\filldraw[black] (3,0) circle (2.5pt) node[anchor=west]{};
\filldraw[black] (1,-3) circle (2.5pt) node[anchor=west]{};
\filldraw[black] (5,-3) circle (2.5pt) node[anchor=west]{};
\end{scope}
\end{tikzpicture}
\subcaption*{Connector (1): ker$(\pi_i)=\langle f_1 \rangle$}
\end{minipage}
\begin{minipage}{.5\textwidth}

\begin{tikzpicture}

\begin{scope}state without output/.append style={draw=none}%[every node/.style={circle,thick,draw}]
    \node (A) at (3,0) { };
    \node (B) at (0,-3) { };

    \node (C) at (6,-3) { };
   
\end{scope}

\begin{scope} %[>={Stealth[black]},
              every node/.style={fill=white,circle},
              every edge/.style={draw=black,very thick}]
              
    %\path [->] (A) edge [bend left=40] node {} (B);
    %\path [->] (B) edge [bend left=40] node {} (A);
    %\draw [->] (B) to  [out=320,in=40,looseness=9] (B);
    %\path [->] (D) edge node {$3$} (C);
    %\path [->] (A) edge node {$3$} (E);
    %\path [->] (D) edge node {$3$} (E);
    %\path [->] (D) edge node {$3$} (F);
    %\path [->] (C) edge node {$5$} (F);
    %\path [->] (E) edge node {$8$} (F); 
    %\path [->] (B) edge[bend right=60] node {$1$} (E); 

    %\draw[]  (A) node[draw=none][midway,above] {$b_1$} (B);
    \path [->] (A) edge [bend right=10] node [midway,above] {$f_2$} (B);
   
    \path [->] (B) edge [bend right=10] node [midway,above] {$f_2$} (C);
    \path [->] (C) edge [bend right=10] node [midway,right] {$f_2$} (A);
   
    %\draw [->] (A) to  [out=90,in=140,looseness=8] edge node [midway,left] {$c_2$} (A);
    \path [->] (A) edge [bend left=10] node [midway,below] {$f_1$} (B);
   
    \path [->] (B) edge [bend left=10] node [midway,above] {$f_1$} (C);
    \path [->] (C) edge [bend left=10] node [midway,left] {$f_1$} (A);
    
    \filldraw[black] (3,0) circle (2.5pt) node{}[anchor=west];
\filldraw[black] (0,-3) circle (2.5pt) node{}[anchor=west];
\filldraw[black] (6,-3) circle (2.5pt) node{}[anchor=west];
\end{scope}

\end{tikzpicture}
\subcaption*{Connector (2): ker$(\pi_i)=\langle f_1^2f_2,f_1f_2^2\rangle$}
\end{minipage}

\begin{minipage}{.5\textwidth}
\begin{tikzpicture}
\begin{scope}state without output/.append style={draw=none}%[every node/.style={circle,thick,draw}]
    \node (A) at (3,0) { };
    \node (B) at (0,-3) { };

    \node (C) at (6,-3) { };
   
\end{scope}

\begin{scope} %[>={Stealth[black]},
              every node/.style={fill=white,circle},
              every edge/.style={draw=black,very thick}]
              
    %\path [->] (A) edge [bend left=40] node {} (B);
    %\path [->] (B) edge [bend left=40] node {} (A);
    %\draw [->] (B) to  [out=320,in=40,looseness=9] (B);
    %\path [->] (D) edge node {$3$} (C);
    %\path [->] (A) edge node {$3$} (E);
    %\path [->] (D) edge node {$3$} (E);
    %\path [->] (D) edge node {$3$} (F);
    %\path [->] (C) edge node {$5$} (F);
    %\path [->] (E) edge node {$8$} (F); 
    %\path [->] (B) edge[bend right=60] node {$1$} (E); 

    %\draw[]  (A) node[draw=none][midway,above] {$b_1$} (B);
    \path [->] (A) edge [bend right=10] node [midway,above] {$f_2$} (B);
   
    \path [->] (B) edge [bend right=10] node [midway,above] {$f_2$} (C);
    \path [->] (C) edge [bend right=10] node [midway,right] {$f_2$} (A);
   
    %\draw [->] (A) to  [out=90,in=140,looseness=8] edge node [midway,left] {$c_2$} (A);
    \path [<-] (A) edge [bend left=10] node [midway,below] {$f_1$} (B);
   
    \path [<-] (B) edge [bend left=10] node [midway,above] {$f_1$} (C);
    \path [<-] (C) edge [bend left=10] node [midway,left] {$f_1$} (A);
    
\filldraw[black] (3,0) circle (2.5pt) node[anchor=west]{};
\filldraw[black] (0,-3) circle (2.5pt) node[anchor=west]{};
\filldraw[black] (6,-3) circle (2.5pt) node[anchor=west]{};
\end{scope}

\end{tikzpicture}
\subcaption*{Connector (3): ker$(\pi_i)=\langle f_1f_2 \rangle$}
\end{minipage}

\begin{minipage}{.5\textwidth}

\begin{tikzpicture}
\begin{scope}state without output/.append style={draw=none}%[every node/.style={circle,thick,draw}]
    \node (A) at (3,0) { };
    \node (B) at (1,-3) { };

    \node (C) at (5,-3) { };
   
\end{scope}

\begin{scope} %[>={Stealth[black]},
              every node/.style={fill=white,circle},
              every edge/.style={draw=black,very thick}]
              
    %\path [->] (A) edge [bend left=40] node {} (B);
    %\path [->] (B) edge [bend left=40] node {} (A);
    %\draw [->] (B) to  [out=320,in=40,looseness=9] (B);
    %\path [->] (D) edge node {$3$} (C);
    %\path [->] (A) edge node {$3$} (E);
    %\path [->] (D) edge node {$3$} (E);
    %\path [->] (D) edge node {$3$} (F);
    %\path [->] (C) edge node {$5$} (F);
    %\path [->] (E) edge node {$8$} (F); 
    %\path [->] (B) edge[bend right=60] node {$1$} (E); 

    %\draw[]  (A) node[draw=none][midway,above] {$b_1$} (B);
    \path [->] (A) edge node [midway,above] {$f_1$} (B);
   
    \path [->] (B) edge node [midway,above] {$f_1$} (C);
    \path [->] (C) edge node [midway,left] {$f_1$} (A);
   
    %\draw [->] (A) to  [out=90,in=140,looseness=8] edge node [midway,left] {$c_2$} (A);
    \path [->] (A) edge [out=30, in=150, looseness=15, "$f_2$"] (A);
     \path [->] (B) edge [out=150, in=250, looseness=15, "$f_2$"] (B);
      \path [->] (C) edge [out=310, in=40, looseness=15, "$f_2$"] (C);
\filldraw[black] (3,0) circle (2.5pt) node[anchor=west]{};
\filldraw[black] (1,-3) circle (2.5pt) node[anchor=west]{};
\filldraw[black] (5,-3) circle (2.5pt) node[anchor=west]{};
\end{scope}
\end{tikzpicture}
\subcaption*{Connector (4): ker$(\pi_i)=\langle f_2 \rangle$}
\end{minipage}
\caption{All shapes of connectors of $G_c$.}
\label{figconc}
\end{figure}

\textbf{Now consider the group $G:=\langle G_{a,b}\bigcup G_c \rangle$}, which is a \textit{fragmentation group of $G_0$}. Note that we need to undo the telescoping to recover the action of $G_c$ on $\Omega(\mathsf{B}_w)$. However, all the notations will be consistent. 

It follows that the points $\mu,\sigma,\lambda,\delta$ are purely non-Hausdorff singularities of $G$. Relabel them respectively by $\xi_0,\xi_1,\xi_2,\xi_3$. Their respective orbital graphs $\Gamma_{\xi_i}$ are, for $i=0, 1, 2$, the inductive limits of embeddings of $\mathcal{T}_{1,3k+i}$ to the end of $\mathcal{T}_{1,3(k+1)+i}$ that contains $\xi_{i,3(k+1)+i}$ (the $(3(k+1)+i)$-th truncation of $\xi_i$), and eventually adding loop at the boundary point of $\mathcal{T}_{\xi_i}$, and, for $i=3$, embedding of $\mathcal{T}_{1,n_k}$ to $\mathcal{T}_{1,n_{k+1}}$ corresponding to $w_{n_k}=w_{n_{k+1}}=3$ along the truncations of $\delta$. The graph of germs $\widetilde{\Gamma}_{\xi_i}$, for $i=0,1,2$, is obtained by connecting $4$ copies of $\Gamma_{\xi_i}$ by the Cayley graph of the four-group $H_i:=\{1,b_i,c_i,d_i\}$ at $\xi_i$. For $i=3$, the graph of germs $\widetilde{\Gamma}_{\xi_3}$ is obtained by take $|G_c|$ $(=9)$ copies of $\Gamma_{\xi_3}$ and connecting them at $\xi_3$ by the Cayley graph of $H_3=G_c$. 

\subsection{Subexponential growth estimate}

The subexponential growth estimate of $G$ is very similar to the one in Subsection \ref{dinftyest}. 

There are four types of connectors corresponding to truncations of $\sigma,\lambda,\mu,\delta$, and labeled respectively by fragmentations of $a,a,b,c$. Denote respectively by $l_{a,n,\xi_0},l_{a,n,\xi_1},l_{b,n,\xi_2}$ and $l_{c,n,\xi_3}$ the gaps of consecutive appearances of the corresponding connectors. For instance, if an $(a,\xi_0)$-connector appears on level $n$, then the next $(a,\xi_0)$-connector appears on level $n+l_{a,n,\xi_0}$. It follows that $l_{a,n,\xi_0},l_{a,n,\xi_1},l_{b,n,\xi_2}\leq 6$ and $l_{c,n,\xi_3}$ are uniformly bounded by Condition \ref{bdd33}. Let $l_{c,n,\xi_3}\leq M$. Since the gaps are uniformly bounded, the tiles on level at least $n+30M$ contains all types of connectors of tiles on levels $n$. Let $L=30M$. 

Let $g\in G$. A \textit{rank $n$ traverse of type $h$} is a walk of $g$ starting and ending at a boundary points of $\mathcal{T}_{j,n}$, not touching the boundary points in between, and this walk can be lifted to a central part of a $\widetilde{\Gamma}_{\xi_i}$ starting in the branch $(Id,\Gamma_{\xi_i})$ and ending in the branch $(h,\Gamma_{\xi_i})$. Denote by $\Theta_{n,h}$ the set of rank $n$ traverse of type $h$ and let $\Theta_n=\bigcup\limits_{h\in H_{i}}\Theta_{n,h}$. Let $T_{n,h}=\#\Theta_{n,h}$, and $T_n=\sum\limits_{h\in H_{i}\backslash\{Id\}}T_{n,h}$. Let $F_h(t)=\sum\limits_{n=0}^{\infty}T_{n,h}t^n$ and $F(t)=\sum\limits_{n=0}^{\infty}T_nt^n$. Then a traverse of $\mathcal{T}_{j,n+L}$ must induce all types of traverses of all $\mathcal{T}_{j,n}$. It follows that 
 \[T_{n+L}\leq \sum_{h\neq Id}T_{n,h},\]
\[\implies t^{-L}(F(t)-(T_0+T_1t+T_2t^2+\ldots +T_{L-1}t^{L-1}))\leq \sum\limits_{h\neq Id}F_h(t)\]
\[\implies F(t)\leq |w| \dfrac{g(t)}{Ct^{-L}-D}\]
where $g(t)>0$ for $t>0$ and $0<D<C$. Hence if $t\in (0,\eta^{-1})$, where $\eta=\sqrt[L]{\frac{D}{C}}$ is the real positive root of the polynomial equation $Cx^L-D=0$, the right-hand side of the above inequality is positive.

Now let us estimate the cardinalities $|\mathcal{T}_{j,n}|$. We have
\[|\mathcal{T}_{1,n}|+|\mathcal{T}_{2,n}|\leq|\mathcal{T}_{1,n+1}|\leq 3|\mathcal{T}_{1,n}|+|\mathcal{T}_{2,n}|\] and \[|\mathcal{T}_{2,n+1}|=|\mathcal{T}_{1,n}|.\]
Hence $\limsup_{n\rightarrow\infty}(\sum_{j=1}^2|\mathcal{T}_{j,n}|)^{1/n}$ exists and is denoted $\beta$. Note that $\beta>1$. It follows that the growth function of $G$, denoted $\gamma_G(R)$ is dominated by $\exp(R^{\alpha})$ for every $\alpha>\frac{\log \beta}{\log \beta-\log\eta}$.

\subsection*{Acknowledgements} I would like to thank Volodymyr Nekrashevych for his patient guidance through the completion of this work.

\def\BState{\State\hskip-\ALG@thistlm}
\makeatother
%%%%%%%%%%%%%%%%%%%%%%%%%%%
\let\oldbibitem\bibitem
\renewcommand{\bibitem}{\setlength{\itemsep}{0pt}\oldbibitem}
%%%%%%%%%%%%%%%%%%%%%%%%%%%%%%%%%%%%%%%%%%%%%%%%%%%%%%%%%%%%%%%
%The bibliography style declared is the IEEE format. If
%you require a different style, see the document
%bibstyles.pdf included in this package. This file,
%hosted by the University of Vienna, shows several
%bibliography styles and examples of in-text citation
%and a references page.
%\newpage
%\let\oldbibitem\bibitem
%\renewcommand{\bibitem}{\setlength{\itemsep}{0pt}\oldbibitem}
%\bibliographystyle{ieeetr}

%\addcontentsline{toc}{section}{REFERENCES}

%\renewcommand{\bibname}{{\normalsize\rm REFERENCES}}

%\bibliography{Growth.bib}
%\printbibliography
%\newpage
% \printbibliography[
% heading=bibintoc,
% title={References}
% ] 
% %Prints the entire bibliography with the titel "Whole bibliography"

% \clearpage

%Filters bibliography
%\printbibliography[heading=subbibintoc,type=article,title={Articles only}]
%\printbibliography[type=book,title={Books only}]

%\printbibliography[keyword={physics},title={Physics-related only}]
%\printbibliography[keyword={latex},title={\LaTeX-related only}]
%\appendix
%\bibliographystyle{alpha}
%\bibliography{growthref}

\begin{thebibliography}{BNZ22}

\bibitem[BNZ22]{BNZ}
Laurent Bartholdi, Volodymyr Nekrashevych, and Tianyi Zheng.
\newblock Growth of groups with linear {Schreier} graphs, 2022.
\newblock arXiv:2205.01792.

\bibitem[Bon07]{Bon}
Ievgen Bondarenko.
\newblock Groups generated by bounded automata and their {Schreier} graphs.
\newblock 2007.
\newblock Ph.D. Dissertation. Texas A\&M University.

\bibitem[Can19]{JC19}
Justin Cantu.
\newblock Periodic groups via orbital graphs.
\newblock 2019.
\newblock Doctoral dissertation, Texas A\&M University.

\bibitem[Fek23]{fek23}
M.~Fekete.
\newblock Über die verteilung der wurzeln bei gewissen algebraischen gleichungen mit ganzzahligen koeffizienten.
\newblock {\em Mathematische Zeitschrift}, 17:228--249, 1923.

\bibitem[Gro81]{Gro81}
Michael Gromov.
\newblock Groups of polynomial growth and expanding maps (with an appendix by jacques tits).
\newblock {\em Publications Mathématiques de l'IHÉS}, 53:53--78, 1981.

\bibitem[Kua24]{kua24a}
Zheng Kuang.
\newblock Growth of groups with incompressible elements, \rom{1}, 2024.
\newblock arXiv:2402.16238.

\bibitem[Mil06]{mil2006}
John Milnor.
\newblock {\em Dynamics in One Complex Variable. Third Edition.}
\newblock Princeton University Press, 2006.

\bibitem[Nek18]{nekrash18}
Volodymyr Nekrashevych.
\newblock {Palindromic subshifts and simple periodic groups of intermediate growth}.
\newblock {\em Annals of Mathematics}, 187(3):667---719, 2018.

\bibitem[Nek22]{nek22}
Volodymyr Nekrashevych.
\newblock {\em Groups and Topological Dynamics}.
\newblock Graduate studies in mathematics. American Mathematical Society, 2022.

\end{thebibliography}

\end{document}